\documentclass[a4paper,12pt,oneside]{article}

\usepackage{graphicx}
\usepackage[all]{xy}
\usepackage[centertags]{amsmath}
\usepackage{latexsym}
\usepackage{xcolor}
\usepackage{amsfonts}
\usepackage{tikz-cd}
\usepackage{pst-node}
\usepackage{amssymb,amsthm}
\usepackage[english]{babel}

\usepackage{newlfont}

\usepackage{indentfirst}
\usepackage[english]{babel}
\usepackage[latin1]{inputenc}
\usepackage{newlfont}

\usepackage{indentfirst}
\usepackage[english]{babel}
\usepackage[latin1]{inputenc}
\usepackage{eufrak}
\usepackage{hyperref}

\newtheorem{remark}{Remark}[section]

\theoremstyle{plain}
\newtheorem{lemma}{Lemma}[section]

\newtheorem{proposition}{Proposition}[section]
\newtheorem{theorem}{Theorem}[section]
\newtheorem{definition}{Definition}[section]
\newtheorem{corollary}{Corollary}[section]

\newtheorem{example}{Example}[section]

\newcommand{\G}{\mathcal{G}}

\newcommand{\beqn}{\begin{eqnarray}}
\newcommand{\eeqn}{\end{eqnarray}}
\newcommand{\beq}{\begin{eqnarray}}
\newcommand{\eeq}{\end{eqnarray}}

\newcommand{\bpro}{\begin{proposition}}

\newcommand{\epro}{\end{proposition}}
\newcommand{\blem}{\begin{lemma}}
\newcommand{\elem}{\end{lemma}}
\newcommand{\bdfn}{\begin{definition}}
\newcommand{\edfn}{\end{definition}}
\newcommand{\bcor}{\begin{corollary}}
\newcommand{\ecor}{\end{corollary}}
\newcommand{\bthm}{\begin{theorem}}
\newcommand{\ethm}{\end{theorem}}
\newcommand{\bex}{\begin{example}}
\newcommand{\eex}{\end{example}}
\newcommand{\brmq}{\begin{remark}}
\newcommand{\ermq}{\end{remark}}
\newcommand{\benum}{\begin{enumerate}}
\newcommand{\eenum}{\end{enumerate}}
\newcommand{\bitem}{\begin{itemize}}
\newcommand{\eitem}{\end{itemize}}

\theoremstyle{plain}

\usepackage[latin1]{inputenc}

\title{On dual quaternions, dual split quaternions and Cartan-Schouten metrics on perfect Lie groups}

\author{ Andr\'e Diatta$^{( 1)}$\footnote{ \footnotesize \noindent (1) Aix-Marseille Univ, CNRS, Centrale Marseille, Institut Fresnel, 13013 Marseille, France.
\newline Email: andre.diatta@fresnel.fr; andrediatta@gmail.com.}
;  Bakary Manga$^{( 2)}$
and  Fatimata Sy$^{( 2)}$\footnote{\footnotesize \noindent(2) D\'epartement de Math\'ematiques et Informatique, Universit\'e Cheikh Anta Diop de Dakar, BP 5005 Dakar-Fann, Dakar, S\'en\'egal. Email: bakary.manga@ucad.edu.sn; syfatima89@gmail.com }   }  

\begin{document}
\maketitle

\begin{abstract}
We discuss Cartan-Schouten metrics (Riemannian or pseudo-Riemannian metrics that are parallel with respect to the Cartan-Schouten canonical connection) on perfect Lie groups.  Applications are foreseen in Information Geometry.  
Throughout this work,  the tangent bundle $TG$ and the cotangent bundle $T^*G$ of a Lie group $G$,  {\color{black} are } always endowed with their Lie group structures induced by the right trivialization.
 We show that $TG$ and $T^*G$ are isomorphic if $G$ possesses a biinvariant Riemannian or pseudo-Riemannian metric. We also show that, if on a perfect Lie group,  there exists a Cartan-Schouten metric, then it  must be biinvariant. We compute all such metrics on the cotangent bundles of simple Lie groups. We further show the following.
Endowed with their canonical Lie group structures,  the set of unit dual quaternions is isomorphic   to $T^*SU(2)$, the set  of unit dual split quaternions is isomorphic to $T^*SL(2,\mathbb R)$. 
 The group SE(3) of special  rigid displacements of the Euclidean  $3$-space is isomorphic to $T^*$SO(3). The group $SE(2,1)$ of special rigid displacements of the Minkowski $3$-space is isomorphic to $T^*$SO(2,1). 
  Some results on SE($3$) by N. Miolane and X. Pennec, and   M. Zefran,  V. Kumar and C. Croke, are generalized  to $SE(2,1)$ and to $T^*G$, for any simple Lie group $G$. 
\end{abstract}

\maketitle
\section{Introduction}

 One of the interests in the canonical Cartan-Schouten symmetric connection on Lie groups,  introduced by E. Cartan  and J. A. Schouten in \cite{Cartan-Schouten}, is  the inverse problem of Lagrangian dynamics for the system of second order differential equations corresponding to the induced geodesic spay (\cite{crampin-mestdag, ghanam-thompson-miller, bi-invariant-and-noninvariant-metrics-ghanam-hindeleh-thompson,muzsnay-thompson,muzsnay,rawashdeh-thompson,thompson2-3D,thompson3D,strugar-thompson}). 
 A particular case  is when there exists at least one Lagrangian which is quadratic in the velocities such that the corresponding Hessian $\mu$ is a metric on the Lie group $G$.  In that case, such a metric $\mu$  is covariantly constant with respect to the connection.

From now on and throughout this work, the word 'metric' refers to both Riemannian and pseudo-Riemannian metrics. 
 In the present paper,  metrics $\mu$  which are covariantly constant $\nabla\mu =0$, where $\nabla$ is the canonical Cartan-Schouten connection, are referred to as Cartan-Schouten metrics (Definition \ref{def:Cartan-Schouten-Metrics}).   Lie groups of dimension $\leq 6$,  possessing a Cartan-Schouten metric have already been classified (\cite{bi-invariant-and-noninvariant-metrics-ghanam-hindeleh-thompson,ghanam-thompson-miller, muzsnay,rawashdeh-thompson,strugar-thompson,thompson2-3D,thompson3D}).
  In \cite{Zefran}, in the framework of kinematic analysis and robot trajectory planning, 
the authors address the problem of finding metrics all of whose geodesics are screw motions 
in the special Euclidean group $SE(3,\mathbb R)$ of rigid motions. Cartan-Schouten metrics turn out to be the only solutions to such a problem.  
On the other hand, a more general statistical framework based on Riemannian Geometry 
has been introduced, and studied by many authors (see e.g. \cite{lauritzen87,Lorentzi-Pennec,Matsuzoe2007,Matsuzoe2010, pennec2}), which takes into account both Euclidean and curved spaces. 
It has been proposed  that such a framework be used on Lie
groups (\cite{gallier-quaintance, miolane-pennec2015, pennec1}). This entails that the underlying metric be compatible with the group structure. 
So far, the researchers and experts have  linked such a compatibility  with the condition that the metric should be biinvariant, that is, invariant under both left and right translations of the Lie group. 
However, as pointed out in \cite{miolane-pennec2015},  some of the Lie groups which are commonly used in some applications do not carry such metrics.
To overcome such an apparent limitation, we propose (see \cite{diatta-manga-sy-jmp, FatimataSy}) to drop the invariance assumption on the metric. We only require that the metric abides by certain Lie group features. Namely, we require that 1-parameter subgroups be geodesics, 
the Riemannian mean coincide  with the Lie group biinvariant exponential barycenter as in \cite{pennec-biinvariant-means}, etc. We thus propose the natural framework of Cartan-Schouten metrics.
From  the geometric viewpoint, Cartan-Schouten metrics are the natural  generalization of biinvariant metrics as in \cite{bi-invariant-and-noninvariant-metrics-ghanam-hindeleh-thompson, medina85, medina-revoy-ENS, pennec-biinvariant-means}. 

 In this paper, we prove the following.
 Endowed with their Lie group structures given by the right trivialization, the tangent bundle $TG$ and the cotangent bundle $T^*G$ of a Lie group $G$, are isomorphic if $G$ itself possesses a biinvariant  metric (Theorem \ref{thm:cotangentbundlerighttrivialization}). The converse will be discussed elsewhere.
We prove that every Cartan-Schouten metric on a perfect Lie group is necessarily biinvariant (Theorem \ref{Theorem:parallel-metric-perfect-Lie-groups}. In particular, every Cartan-Schouten metric on the cotangent bundle of a simple Lie group is necessarily biinvariant. We further compute all such biinvariant metrics on $T^*G$ when G is a simple Lie group (Theorem \ref{thm:biinvariant-metric-on-dual}). 
 With their canonical Lie group structures,  the set of unit dual quaternions  (resp. unit dual split quaternions) is isomorphic to $T^*SU(2)$ (resp. $T^*SL(2,\mathbb R)$). The group SE(3) of special  rigid displacements of the Euclidean  $3$-space is isomorphic to  the Lie group $T^*$SO(3). The group $SE(2,1)$ of rigid displacements of the Minkowski $3$-space is isomorphic to the Lie group $T^*$SO(2,1).  See Theorems \ref{theorem:commutingdiagramSO(3)}, \ref{theorem:cummutingdiagramSL(2)}.
In the mathematical setting, this work also implies, in particular, that the sets of unit dual quaternions, unit dual split quaternions and of the rigid displacements  $SO(3)\ltimes \mathbb R^3$, $SO(2,1)\ltimes \mathbb R^3$,  
are all symplectic manifolds and as Lie groups, they all carry biinvariant pseudo-Riemannian structures.  
The isomorphism between $T^*SO(3)$ and $SE(3)$ could be used to look at the result on Cartan-Schouten metrics on cotangent bundles of simple Lie groups, as a generalization of some results on SE($3$) by 
N. Miolane and Pennec (\cite{miolane-pennec2015}, \cite{pennec1}), and   M. Zefran,  V. Kumar and C. Croke (\cite{Zefran}). 
The equivalences proved here also open more routes towards applications of dual quaternions and dual split quaternions in many areas such as thermodynamics, statistical mechanics, information geometry and machine learning. 

The paper is organized as follows. 
Section \ref{chap:Onstatistical-structures-on-groups} is devoted to a few definitions, examples, a characterization of Lie groups with a biinvariant metrics (Subsection \ref{chap:Liegroups-with-biinvariant-metric}) and Cartan-Schouten metrics on perfect Lie groups (Subsection \ref{chap:Cartan-Schouten-metrics-on-perfect-Liegroups}).  
We discuss the Cartan-Schouten metrics on cotangent bundles of simple Lie groups in Section \ref{chap:Cartan-Schouten-of-cotangent-bundles-of-simple-Lie-groups}.
Section \ref{chap:reminders-quaternions-plit-quaternions} is a reminder on quaternions, dual quaternions, split quaternions, dual split quaternions and some of their relationships with the groups of rigid motions.
In Section \ref{chap:rigid-motions-dual}, we discuss the isomorphisms between the cotangent bundles and the groups of rigid motions, as well as the problem of metrics all of whose geodesics are screw motions.
 
We will let $E_{i,j}$ stand for the elementary $n\times n$ matrix with zero in all entries except the $(i,j)$ entry which is equal to $1$.
If $V$ is a vector space, we denote its  linear dual by $V^*$ and by $\mathbb I_{V}$ its identity map.  We let $(e_1^*,\dots,e_n^*)$ stand for the dual basis of a basis $(e_1,\dots,e_n)$. All Lie groups and Lie algebras considered here are real.

\section{On Cartan-Schouten metrics on Lie groups}\label{chap:Onstatistical-structures-on-groups}

Let us remind that Cartan-Schouten connections on a Lie group $G$, are the left invariant connections whose geodesics through the identity are $1$-parameter subgroups of $G$. 
The most prolific examples are the classical  $+$, $-$ and $0$ Cartan-Schouten connections (\cite{Cartan-Schouten})  respectively given in the Lie algebra $\mathcal G$ of $G$ by
\beqn\label{Cartan-connections} \nabla_x y:=\lambda [x,y], 
\eeqn
$ \lambda=1,0,\frac{1}{2}$, for all  $x,y\in\mathcal G$.
\begin{definition} 
The $0$-connection given by  $\nabla_xy=\frac{1}{2} [x,y]$,  for all $x,y\in\mathcal G$,  is also termed the Cartan-Schouten canonical connection. It is the unique symmetric (torsion free) Cartan-Schouten connection which is  bi-invariant.
\end{definition}
Consider a metric $\mu$ on $G$,
 which is covariantly constant (or equivalently, parallel) with respect to the Cartan-Schouten 
canonical connection $\nabla$. 
The latter property simply reads $\nabla \mu=0$  or, equivalently,
\beqn \label{eq:parallel} 
x^+\cdot\mu(y^+,z^+)=\frac{1}{2}\Big(\mu([x^+,y^+],z^+)+\mu(y^+,[x^+,z^+])\Big) 
\eeqn
for any left invariant vector fields $x^+,y^+,z^+$ on $G$. 
\begin{definition} \label{def:Cartan-Schouten-Metrics} If a metric  on a Lie Group, is parallel with respect to the Cartan-Schouten canonical connection, we call it a Cartan-Schouten metric.
\end{definition}

\subsection{Example: The Heisenberg  Lie group $\mathbb H_3$ of dimension 3}\label{Heisenberg}
Consider the  
 Heisenberg group  $\mathbb H_3:=\left\{
\begin{pmatrix}1&x&z
\\ 0&1&y\\ 
0&0&1
\end{pmatrix}, \; x, y, z\in\mathbb R\right\}$. Its Lie algebra $\mathcal{H}_3$ is spanned by the $3\times 3$ elementary matrices  $e_1:=E_{1,2}$, $e_2:=E_{2,3}$, $e_3:=E_{1,3}$.  So the Lie bracket reads $[e_1,e_2]=e_3$. We identify $\mathbb H_3$ with $\mathbb R^3$, with the multiplication 
 $(x,y,z) (x',y',z')=(x+x',y+y',z+z'+xy')$.
According to  \cite{medina85}, \cite{medina-revoy-ENS}, the existence of a biinvariant metric on G, with Lie algebra $\mathcal G$ and center $Z(\mathcal G),$ implies the equality $\dim [\mathcal G,\mathcal G]+\dim Z(\mathcal G)=\dim \mathcal G.$ 
So, $\mathbb H_3$ does not have any biinvariant metric, since $[\mathcal H_3, \mathcal H_3]=Z(\mathcal H_3)=\mathbb R e_3.$ 
However,  $\mathbb H_3$ does possess infinitely many  Cartan-Schouten metrics. 
\begin{proposition}
 Any Cartan-Schouten  metric $\mu$ on $\mathbb H_3,$  is of the form
\begin{eqnarray}
 \mu\!\!\!&=&\!\!\left(\frac{1}{4}ay^2-cy+m\right)dx^2 +\left(\frac{1}{4}ax^2-bx+e\right)dy^2 +a\,dz^2  \cr 
& &\!\!\! +\!\left(\frac{1}{4}axy\!-\!\frac{1}{2}cx\!-\!\frac{1}{2}by\!+\!d\right)\, dxdy
\!-\! \left(\frac{1}{2}ay\!-\!c\right)\, dxdz \!-\! 
\left(\frac{1}{2}ax\!-\!b\right) \, dydz,\nonumber
\end{eqnarray}
with  $aem-ad^2-b^2m+2bcd-c^2e\neq 0,$ where $a,b,c,d,e,m$ are real constants.
\end{proposition}
\brmq 
(A) If we set $a=e=m=1$ and $b=c=d=0$, we recover the metric given in 
\cite{bi-invariant-and-noninvariant-metrics-ghanam-hindeleh-thompson} : 
$ 
\mu=dx^2 +dy^2 +\left(dz-\frac{y}{2}dx-\frac{x}{2}dy\right)^2,
$ which is a Riemannian metric. 
(B) For $m=e=-a=1$ and $b=c=d=0$ we get
$ 
\mu=( 1-\frac{1}{4}y^2)dx^2 +( 1-\frac{1}{4}x^2)\,dy^2-dz^2 - \frac{1}{4}xy\,dxdy + \frac{1}{2}y\,dxdz+  \frac{1}{2}x\,dydz$, which is a Lorentzian metric. Note that $\det(\mu)=-1$ and $\mu(\vec{v},\vec{v})=0$, where $\vec{v}=\frac{\partial}{\partial y}+(\frac{x}{2}+1)\frac{\partial}{\partial z}$. 
Both metrics in (A) and in (B) have the same Levi-Civita connection, although one is Riemannian and the other Lorentzian.
\ermq

\subsection{Example:  Lorentz Lie group SO(3,1) and Lie algebra $\mathfrak{so}(3,1)$}\label{LorentzGroup}

Let $G$ be the group of Lorentz transformations (isometries of space-time) of Minkowski 4-space, with determinant $+1$. 
A basis of its Lie algebra $\mathcal G$ is made of the $4\times 4$ matrices $S_1=E_{1,4}+E_{4,1}$,   $S_2=E_{2,4}+E_{4,2}$,  $S_3=E_{3,4}+E_{4,3}$, 
 $S_4=E_{2,3}-E_{3,2}$,  $S_5=E_{3,1}-E_{1,3}$,  $S_6=E_{2,1}-E_{1,2}$.  
 Its Lie bracket reads
$[S_1,S_2]= - S_6$, $[S_1,S_3]= - S_5$, $[S_1,S_5]= - S_3$, $[S_1,S_6]= - S_2$,
 $[S_2,S_3]= S_4$, $[S_2,S_4]= S_3$, $[S_2,S_6]=S_1$,
 $[S_3,S_4]=-S_2$, $[S_3,S_5]=S_1$, 
$[S_4,S_5]=S_6$, $[S_4,S_6]= - S_5$, $[S_5,S_6]=S_4$.
Since $G$ is a simple Lie group, Theorem \ref{Theorem:parallel-metric-perfect-Lie-groups} ensures that every Cartan-Schouten metric  $\mu$ on $G$ is biinvariant. A direct calculation shows that the matrix of $\mu$, in the basis 
$(S_j)$, is 
\beqn (\mu(S_i,S_j))&=&k_1(E_{1,1}+E_{2,2}+E_{3,3}-E_{4,4}-E_{5,5}-E_{6,6})\nonumber\\
&& +  k_2(E_{1,4}+E_{4,1}+E_{2,5}+E_{5,2}-E_{3,6}-E_{6,3})\eeqn
$k_1,k_2\in\mathbb R$. Note that  
$4(E_{1,1}+E_{2,2}+E_{3,3}-E_{4,4}-E_{5,5}-E_{6,6})$ is the matrix of the Killing form $K_0$ and $E_{1,4}+E_{4,1}+E_{2,5}+E_{5,2}-E_{3,6}-E_{6,3}$ is the matrix of $K_J:=K_0(J(\cdot),\cdot)$
where $J=E_{1,4}+E_{4,1}+E_{2,5}+E_{5,2}-E_{3,6}-E_{6,3}$ satisfies $J^2=-\mathbb I_{\mathcal G}$ and $J[x,y]=[Jx,y]$, for any $x,y\in\mathcal G$.

\subsection{A characterization of Lie groups with a biinvariant metric}\label{chap:Liegroups-with-biinvariant-metric}
Recall that  the trivialization by right translations, or simply the right trivialization, of the cotangent bundle  $T^*G$ of a Lie group $G$, is given by the isomorphism $\zeta$ of vector bundles,
$\zeta: T^*G \to G\times \G^*$, $(\sigma,\nu_\sigma)\mapsto (\sigma,\nu_\sigma\circ T_\epsilon R_\sigma)$, 
 where $R_\sigma$ is the right multiplication  $R_\sigma(\tau):=\tau \sigma$ by $\sigma$ in $G$ and $T_\epsilon R_\sigma$ is the derivative of $R_\sigma$ at the unit $\epsilon$ of $G$. 
Likewise, the right trivialization of the tangent bundle $TG$ is given by the
isomorphism $\xi$ of vector bundles,
$\xi: TG \to G\times \G$, $(\sigma,X_\sigma)\mapsto (\sigma, T_\sigma R_{\sigma^{-1}} X_\sigma)$.
Let $G\ltimes_{Ad} \mathcal G$ and  $G\ltimes_{Ad^*} \mathcal G^*$  (or simply $G\ltimes \mathcal G$ and  $G\ltimes \mathcal G^*$, for short) stand for the manifolds $G\times \mathcal G$ and  $G\times \mathcal G^*$, endowed with their Lie group structures  respectively defined by the following products
\begin{eqnarray}
(\sigma_1,x)(\sigma_2,y)&:=&\Big(\sigma_1\sigma_2, x+ Ad_{\sigma_1} y\Big), \label{Adjoint0} \\ 
 (\sigma_1,f)(\sigma_2,g)&:=&\Big(\sigma_1\sigma_2,f+Ad_{\sigma_1}^* g\Big); \label{Adjoint}
 \end{eqnarray}
for any $(\sigma_1,x), (\sigma_2,y)\in G\times \mathcal G$ and any 
$(\sigma_1,f), (\sigma_2,g)\in G\times \mathcal G^*$. 
Here,  $Ad:G\times \mathcal G\to \mathcal G$  and  $Ad^*:G\times \mathcal G^*\to \mathcal G^*$ are respectively the adjoint and coadjoint action of  $G$ on $\mathcal G$ and $\mathcal G^*$, given by $Ad_{\sigma}x=\frac{d}{dt}\vert_{t=0}\sigma \exp(tx) \sigma^{-1}$ and 
$Ad_{\sigma}^*f:=f\circ Ad_{\sigma^{-1}}$,  for $\sigma\in G $, $x\in\mathcal G$, $f\in\mathcal G^*$  and $\exp: \mathcal G \to G$ is the Lie group exponential.  
Thus, $T^*G$ inherits a Lie group structure obtained by pulling back the product (\ref{Adjoint}) in $G\ltimes \mathcal G^*$ via $\zeta$, so that 
\begin{eqnarray}\label{coadjoint2}(\sigma,\nu_\sigma)(\tau,\alpha_\tau)&:=&\zeta^{-1}\Big(\zeta(\sigma,\nu_\sigma)\zeta(\tau,\alpha)\Big)
\nonumber\\
& =&\Big( \sigma\tau,\nu_\sigma\circ  T_{\sigma\tau} R_{\tau^{-1}}\;+\;\; \alpha_\tau\circ T_{\sigma\tau} L_{\sigma^{-1}}\Big).
\end{eqnarray}
Likewise, $TG$ also inherits a Lie group structure, the pullback of (\ref{Adjoint0}) using $\xi$, so that 
\begin{eqnarray}\label{adjoint1}(\sigma,X_\sigma)(\tau,Y_\tau):=\xi^{-1}\Big(\xi(\sigma,X_\sigma)\xi(\tau,Y_\tau)\Big) = \Big( \sigma\tau, T_\sigma R_{\tau}X_\sigma +  T_\tau L_{\sigma} Y_\tau\Big).
\end{eqnarray}
\begin{theorem}\label{thm:cotangentbundlerighttrivialization}
Let $G$ be a Lie group,  $T^*G$ and $TG$  its cotangent and tangent bundles endowed with their respective Lie group structures (\ref{coadjoint2}) and  (\ref{adjoint1}). 
 If $G$ has a biinvariant metric, then the Lie groups $T^*G$ and $TG$ are isomorphic. 
\end{theorem}
\begin{proof}
Let $G$ be a Lie group endowed with a biinvariant metric, say $\mu$.  By abuse of notation, we denote again by $\mu$ its value $\mu_\epsilon$ at the unit  $\epsilon$ of $G$. Define the linear invertible map $\Theta:\mathcal G \to \mathcal G^*, $ $x\mapsto \Theta(x)$, where $\Theta(x):=\mu(x,\cdot)$ is the linear form that maps every $y\in\mathcal G$ to $\langle \Theta(x),y\rangle = \mu(x,y)$.  
 The property $\mu(Ad_\sigma x, Ad_\sigma y) =\mu(x,y)$ is equivalent to 
 $\Theta(Ad_\sigma x) = Ad_\sigma^*\Theta(x)$,  for every $\sigma\in G$ and $x,y\in\mathcal G$.
 Let $\Phi: TG\to T^*G$ be the  isomorphism of vector bundles defined by 
$$\Phi(\sigma,X_\sigma)= \Big(\sigma, \Theta(T_{\sigma}R_{\sigma^{-1}}X_{\sigma})\circ T_{\sigma}R_{\sigma^{-1}}\Big) =  \Big(\sigma, \Theta(T_{\sigma}L_{\sigma^{-1}}X_{\sigma})\circ T_{\sigma}L_{\sigma^{-1}}\Big). $$
 From the definition of $\Phi$, the image
$\Phi\Big((\sigma,X_\sigma)(\tau,Y_\tau)\Big)$ of (\ref{adjoint1}) comes to 
\beqn
\Phi\Big((\sigma,X_\sigma)(\tau,Y_\tau)\Big)&=&
\Big( \sigma\tau\; , \; \Theta( T_{\sigma\tau} R_{\tau^{-1} \sigma^{-1}} \; T_\sigma R_{\tau}X_\sigma )\circ T_{\sigma\tau} R_{\tau^{-1} \sigma^{-1}}
\nonumber\\
 & & \: \qquad +\:
 \Theta( T_{\sigma\tau} R_{\tau^{-1} \sigma^{-1}} \; T_\tau L_{\sigma} Y_\tau)\circ T_{\sigma\tau} R_{\tau^{-1} \sigma^{-1}}\Big). \nonumber
\eeqn 
We further use the equalities 
$T_{\sigma\tau} R_{\tau^{-1} \sigma^{-1}} \; T_\sigma R_{\tau}=T_{\sigma} R_{\sigma^{-1}}$ and 
\beqn
\Theta( T_{\sigma\tau} R_{\tau^{-1} \sigma^{-1}} \; T_\tau L_{\sigma} Y_\tau)\circ T_{\sigma\tau} R_{\tau^{-1} \sigma^{-1}}=
\Theta( T_{\tau} L_{\tau^{-1}} Y_\tau)\circ T_{\sigma\tau} L_{\tau^{-1} \sigma^{-1}}.
\eeqn 
So now $\Phi\Big((\sigma,X_\sigma)(\tau,Y_\tau)\Big)$ becomes
\beqn 
\Big( \sigma\tau\;,\;  \Theta( T_{\sigma} R_{\sigma^{-1}} X_\sigma )\circ T_{\sigma\tau} R_{\tau^{-1} \sigma^{-1}}
\;+\;
 \Theta( T_{\tau} L_{\tau^{-1}} Y_\tau)\circ T_{\sigma\tau} L_{\tau^{-1} \sigma^{-1}}\Big). \nonumber
 \eeqn

\noindent On the other hand, the product $\Phi(\sigma,X_\sigma)\Phi(\tau,Y_\tau)$ is equal to
\beqn
  \Big( \sigma\tau\, ,\, \Theta(T_{\sigma}R_{\sigma^{-1}}X_{\sigma})\!\circ\! T_{\sigma}R_{\sigma^{-1}}\!\circ\! 
  T_{\sigma\tau} R_{\tau^{-1}}\!
+\! \;\Theta(T_{\tau}L_{\tau^{-1}}Y_{\tau})\!\circ\! T_{\tau}L_{\tau^{-1}} \!\circ\! T_{\sigma\tau} L_{\sigma^{-1}} \Big)\nonumber
\eeqn 
\noindent 
and visibly, the latter coincides with $\Phi\Big((\sigma,X_\sigma)(\tau,Y_\tau)\Big)$.
\end{proof} 
The Lie algebra $T^*\mathcal G$ of $T^*G$, is then isomorphic to the semi-direct {\color{black}sum} $\G\ltimes\G^*$, the Lie bracket of two elements $(x,f)$ and $(y,g)$ of $T^*\mathcal G$ being
\beq \label{cotangent-bundle-Lie-bracket}
[(x,f),(y,g)]:=([x,y],ad^*_xg-ad^*_yf).
\eeq
Similarly, the Lie algebra $T\mathcal G$ of $TG$, is the semi-direct sum $\G\ltimes\G$, the Lie bracket of two elements $(x_1,y_1)$ and $(x_2,y_2)$ being 
\beqn[(x_1,y_1),(x_2,y_2)]=([x_1,x_2], [x_1,y_2]-[x_2,y_1]).
\eeqn

\subsection{Cartan-Schouten metrics on perfect Lie groups}\label{chap:Cartan-Schouten-metrics-on-perfect-Liegroups}

In this section, we prove the following (without loss of generality, we can work with connected Lie groups).

\begin{theorem} \label{Theorem:parallel-metric-perfect-Lie-groups} Let $G$ be a perfect Lie group, that is, its Lie algebra $\mathcal G$  satisfies $[\mathcal G,\mathcal G]=\mathcal G$.  
If $G$ possesses a Cartan-Schouten metric $\mu$, then $\mu$ is necessarily a biinvariant metric.
 \end{theorem}
\begin{proof}From Ambrose-Singer holonomy theorem, the curvature tensor $R^\nabla$ {\color{black} of $\mu$} must be skew-symmetric with respect to $\mu$. That is,
\beqn\label{skew-symmetric-curvature}
\mu\Big(R^\nabla (X_1,X_2)Y,Z\Big)+\mu\Big(Y,R^\nabla (X_1,X_2)Z \Big) =0,
\eeqn
for any smooth vector fields $X_1,X_2,Y,Z$ on $G$.  In particular,  if $X_1,X_2,Y,Z$ are all left invariant vector fields on $G$, Equality (\ref{skew-symmetric-curvature}) becomes 
\beqn\label{skew-symmetric-curvature-invariant}
\mu\Big([[X_1,X_2],Y],Z\Big)+\mu\Big(Y,[[X_1,X_2],Z] \Big) =0.
\eeqn
Since $[\mathcal G,\mathcal G]=\mathcal G$, the equality (\ref{skew-symmetric-curvature-invariant}) linearly extends to 
\beqn\label{skew-symmetric-curvature-invariant2}
\mu\Big([X,Y],Z\Big)+\mu\Big(Y,[X,Z] \Big) =0,
\eeqn
for any left invariant vector fields $X,Y,Z$ on $G$. According to (\ref{eq:parallel}), the left hand-side of (\ref{skew-symmetric-curvature-invariant2}) is precisely $2 X\cdot \mu(Y,Z)$. So, the differential of  the function $\mu(Y,Z)$ vanishes, for any left invariant vector fields $Y,Z$ on $G$. Equivalently,  $\mu$ is constant on left invariant vector fields and hence $\mu$ is itself left invariant on $G$, in addition to satisfying (\ref{skew-symmetric-curvature-invariant2}). Consequently, $\mu$ is biinvariant.
\end{proof}

\section{Cartan-Schouten metrics on cotangent bundles of simple Lie groups}\label{chap:Cartan-Schouten-of-cotangent-bundles-of-simple-Lie-groups}

In what follows, if $G$ is a Lie group,  $\mathcal G$ its Lie algebra, the Killing form $K_0$ of $\mathcal G$ will also be looked at as a bilinear  form on the Lie algebra $T^*\mathcal G$ of $T^*G$, 
 with $K_0\Big((x,f),(y,g)\Big)=K_0(x,y)$ for any $(x,f),(y,g)\in T^*\mathcal G$.
The following Lemma (see \cite{medina85}, for part (B)) will be useful in the proof of Theorem \ref{thm:biinvariant-metric-on-dual}. 
Example \ref{LorentzGroup}  is a typical illustration of part (B)-(2).

\begin{lemma}\label{ad-invariant-endomorphisms}Let $\mathcal G$ be a simple  real Lie algebra, $K_0$ its Killing form and  $\dim \mathcal G=n$.
Denote by $K(\mathcal G)$ the space of  linear maps $L:\mathcal G\to \mathcal G$ satisfying $L[x,y]=[L(x),y]$, for any $x,y\in\mathcal G$. 
(A) Every  element of $K(\mathcal G)$ is symmetric with respect to $K_0$. 
(B)   $\dim_{\mathbb R} K(\mathcal G)\leq 2$ and $K(\mathcal G)$ is a commutative subfield of the ring End($\mathcal G$) of endomorphims of the (real) vector space underlying $\mathcal G$.
(1) If $n$ is odd, then $K(\mathcal G)=\mathbb R\mathbb I_{\mathcal G}$. In particular, up to a constant factor,  $K_0$ is the unique ad-invariant metric on $\mathcal G$. 
(2) If $n$ is even, then $\dim K(\mathcal G)=2$. More precisely, $K(\mathcal G)$ is isomorphic to $\mathbb C$, it is spanned by    $\mathbb I_{\mathcal G}$ and some $J$ with $J^2=-\mathbb I_{\mathcal G}$.
\end{lemma}

\begin{proof} (A) Write $L=L_1+L_2$, where
\beqn
 K_0(  L_{1}(x),y) &=&\frac{1}{2} (K_0( L(x),y)+ K_0( x, L(y))), \text{  and }
\nonumber\\
 K_0(  L_{2}(x),y) &=&\frac{1}{2} (K_0( L(x),y)- K_0( x, L(y))),
\nonumber
\eeqn
 for any $x,y\in\mathcal G$. The following hold $  L_{2}[x,y]= [ L_{2}(x),y]$ and 
\beqn
 K_0(  L_{2}([x,y]),z) &=&  K_0(   L_{2}(x),[y,z])  = -  K_0(  x,[ L_{2}(y),z]) \nonumber\\
&=&   -  K_0(  [x, L_{2}(y)],z) = -  K_0(  L_{2}( [x,y]),z),
\eeqn
 for any $x,y,z\in\mathcal G$. In other words $ L_{2}( [x,y])=0$, for any $x,y\in\mathcal G$. Hence, we have $ L_{2}=0$ and $L=L_1 $ is thus symmetric with respect to $K_0$. 
Part (B) has been proved in \cite{medina85}, Sect. 5, p. 409-410.
 \end{proof}
\blem \label{lem:perfectT*G} Let $G$ be a perfect Lie group, $\mathcal G$ its Lie algebra. Suppose $G$ has a Cartan-Schouten metric. Then  $T^*G$ is a perfect Lie group. In particular, if $G$ is semisimple, then $T^*G$ is a perfect Lie group.
\elem
\begin{proof} Let $\mu$ be a Cartan-Schouten metric on $G$. From Theorem \ref{Theorem:parallel-metric-perfect-Lie-groups},  $\mu$ can be seen as an ad-invariant metric on $\mathcal G$.
 The invertible linear map
 $\hat \Theta: \mathcal G\to\mathcal G^*$, $x\mapsto\mu(x,\cdot)$, satisfies $\hat \Theta ([x_1,x_2])=ad_{x_1}^*\hat \Theta (x_2)$, and as $[\mathcal G,\mathcal G]=\mathcal G$, thus for any $\bar g\in\mathcal G^*$ we can suppose, without any loss of generality, that there exist $x_1,x_2\in \mathcal G$ such that 
$\bar g=\hat \Theta ([x_1,x_2])=ad_{x_1}^*\hat \Theta (x_2)=ad_{x}^* g=[x,g]$, where $x_1=: x$ and $\hat \Theta (x_2)=:g$. So in particular, $[\mathcal G,\mathcal G^*]=\mathcal G^*$  and the Lie algebra  $T^*\mathcal G$ of $T^*G$ satisfies  $[T^*\mathcal G,T^*\mathcal G]=T^*\mathcal G$.
\end{proof}
 For an integer $p\ge 0$, let  $\mathbb I_{p,n}$ stand for the diagonal  $n\times n$ matrix with $-1$ in the first $p$ diagonal entries and  $1$ in the last $n-p$  entries.  As above, for  a Lie algebra $\mathcal G$, we let  $K_0$ stand for its Killing form and $\langle,\rangle$ the duality paring between  $\mathcal G$ and $\mathcal G^*$. 
For a linear map $J:\mathcal G\to \mathcal G$  as in (B)-(2) in Lemma \ref{ad-invariant-endomorphisms}, we define the symmetric bilinear forms  $K_J$ and $\langle,\rangle_J$  on $T^*\mathcal G$ as  $K_J(x,y)=K_0(J x,y)$, $K_J(x,g)=K_J(f,g)=0$,
 $\langle x,y\rangle_J =\langle f,g\rangle_J =0$ and $\langle x,g\rangle_J =\langle J(x),g\rangle$, for $x,y\in \mathcal G$, $f,g\in\mathcal G^*$. 
\begin{theorem}\label{thm:biinvariant-metric-on-dual}
Let $G$ be an $n$-dimensional simple Lie group,  $\mathcal G$ its Lie algebra. 
Then every Cartan-Schouten metric $\mu$ on $T^*\mathcal G$ is biinvariant and has signature $(n,n)$.
\bitem 
\item[(A)] If $n$ is odd, then $\mu$ is a linear combination $\mu=sK_0+t\langle,\rangle$ of  $K_0$ and $\langle,\rangle$, with $s,t\in\mathbb R$, $t\neq 0$.
Let $(p,n-p)$ be the signature of $K_0$ and $(e_1,\dots,e_n)$ a basis of $\mathcal G$ in which the matrix of $K_0$ is 
 $\mathbb I_{p,n}$.  
Then $\mu$ has a matrix of the form
$\begin{pmatrix} s\mathbb I_{p,n}&t \mathbb I_{n}\\ t \mathbb I_{n}&{\mathbf 0}_{n}\end{pmatrix}$  in 
the basis $(e_1,\dots,e_n, e_1^*,\dots,e_n^*)$ of  $T^*\mathcal G$, where $s,t\in\mathbb R$, $t\neq 0$ and  ${\mathbf 0}_{n}$ is
 the zero $n\times n$ matrix.
\item[(B)] If $n$ is even, then the space of Cartan-Schouten metrics on $T^*G$ is 4-dimensional. More precisely, each such metric is a linear combination $\mu=s_1K_0+s_2K_J+t_1\langle ,\rangle+t_2\langle,\rangle_J$, with $s_1,s_2,t_1,t_2\in\mathbb R$, with $J$ as in Lemma \ref{ad-invariant-endomorphisms}. 
That is, for any $x,y\in\mathcal G$,  $f,g\in\mathcal G^*$,  
\eitem 
\beqn 
\mu\Big((x,f),(y,g)\Big)&=& s_1K_0(x,y)+s_2K_J(x,y) + t_1\langle (x,f),(y,g)\rangle \cr 
& & +\:  t_2\langle (x,f),(y,g)\rangle_J.
\eeqn
\end{theorem}
\proof  
Lemma \ref{lem:perfectT*G}, ensures that $T^*G$  is a perfect Lie group. {\color{black} Thus,  according to Theorem \ref{Theorem:parallel-metric-perfect-Lie-groups},  any Cartan-Schouten metric $\mu$ on $T^* G,$  must be  biinvariant. } So the restriction of $\mu$ to $T^*\mathcal G$ is ad-invariant.
On $T^*\mathcal G=\mathcal G\ltimes \mathcal G^*$ with its Lie bracket (\ref{cotangent-bundle-Lie-bracket}), an ad-invariant metric $\mu$ is given by a linear invertible map $\phi: T^*\mathcal G\to T^*\mathcal G$ such that $\mu((x,f),(y,g))= \langle \phi(x,f),(y,g) \rangle=\langle(x,f), \phi(y,g) \rangle$. We write $\phi(x,f)= (\phi_{1,1}(x)+\phi_{2,1}(f), \phi_{1,2}(x)+\phi_{2,2}(f))$ where  $\phi_{1,1}:\mathcal G\to\mathcal G$,  $\phi_{1,2}:\mathcal G\to\mathcal G^*$, $\phi_{2,1}:\mathcal G^*\to\mathcal G$ and $\phi_{2,2}:\mathcal G^*\to\mathcal G^*$, are linear maps. The ad-invariance of $\mu$ is equivalent to the relations 
\begin{eqnarray}\label{eq:ad-invariance}
\phi([(x,f),(y,g)]) = [\phi(x,f),(y,g)]= [(x,f),\phi(y,g)],
\end{eqnarray} 
true for any $(x,f)$ and $(y,g)$ in $T^*\mathcal G$. 
The latter equations,  in the case $f=g=0$,  become  
$\phi[x,y] = (\phi_{1,1} [x,y], \phi_{1,2}([x,y]))=
[\phi(x),y] = [\phi_{1,1}(x),y] -  ad_{y}^* \phi_{1,2}(x)$. 
  A component-wise comparison leads to the equalities $ \phi_{1,1} [x,y]=  [\phi_{1,1}(x),y] 
$ and $\phi_{1,2}([x,y])=-  ad_{y}^* \phi_{1,2}(x)
$.  The bilinear form
 $\mu_{\phi_{1,2}}$ on $\mathcal G$  given by 
 $\mu_{\phi_{1,2}} (x,y):=\langle \phi_{1,2}(x),y \rangle$, is ad-invariant. Thus it is of the form $\mu_{\phi_{1,2}} (x,y)=K_0(\psi(x),y)$ for some linear map $\psi$ satisfying $\psi([x,y])=[\psi(x),y]$, for every $x,y\in\mathcal G$.
Now taking $y=0,f=0$ in (\ref{eq:ad-invariance}), leads to 
$\phi[x,g] = \phi_{2,1} (ad_x^*g) + \phi_{2,2}(ad_x^*g)=[\phi(x),g]= ad_{\phi_{1,1}(x)}^*g$,
 or equivalently, $\phi_{2,1} (ad_x^*g) =0$ and $\phi_{2,2}(ad_x^*g)= ad_{\phi_{1,1}(x)}^*g$,
 for any $x\in\mathcal G$ and $g\in\mathcal G^*$. 
Since every element of $\mathcal G^*$ can be taken in the form $ad_{x}^* g$, for some $x\in\mathcal G$ 
and $g\in\mathcal G^*$,  this implies the equality $ \phi_{2,1}=0$. So $\mathcal G^*$ is totally isotropic for the metric $\mu$.  As a consequence, $\mu$ is of signature $(n,n)$.   
We also check that $\phi_{2,2}=\phi_{1,1}^T$. Indeed $\langle \phi_{1,1}^Tg,x\rangle = \langle \phi_{1,1}(x), g\rangle = \langle \phi_{1,1}(x)+\phi_{1,2}(x), g\rangle = 
 \langle \phi(x), g\rangle =\mu(x,g)= \langle x, \phi(g)\rangle= \langle x, \phi_{2,2}(g)\rangle$, for 	any $x\in\mathcal G$, $g\in\mathcal G^*$. Now we discuss the above in the two different scenari where $n$ is odd or even.

(A) Suppose $n$ is odd.  By Lemma \ref{ad-invariant-endomorphisms}, we must have $\phi_{1,1} = t\mathbb I_{\mathcal G}$  and $\psi = s\mathbb I_{\mathcal G}$ for some $t,s\in\mathbb R$,   and  thus $\mu_{\phi_{1,2}} = s K_0$, $s\in \mathbb R$
 and $\phi_{2,2}(ad_x^*g)=t ad_{x}^* g$, for any $x\in \mathcal G$ and $g\in\mathcal G^*$,  or equivalently, 
$\phi_{2,2}=t\mathbb I_{\mathcal G^*}$, where $t\in\mathbb R$ and $\mathbb I_{\mathcal G}$, $\mathbb I_{\mathcal G^*}$ are the identity maps of $\mathcal G$ and  $\mathcal G^*$, respectively.  
Altogether, we have 
 $\mu((x,f), (y,g))=  \langle \phi_{1,2}(x),y\rangle+ \langle \phi_{1,1}(x),g\rangle +\langle \phi_{2,2}(f),y\rangle$
 which also reads $\mu((x,f), (y,g))=\mu_{\phi_{1,2}}(x ,y) + t( \langle x,g\rangle +\langle f,y\rangle)
= sK_0(x,y) +  t \langle (x,f),(y,g)\rangle $.
Now choosing a basis $(e_1,\dots,e_n)$ of $\mathcal G$ in which the matrix $K_0(e_i,e_j)$  of $K_0$ is of the form $\mathbb I_{p,n}$,  leads to the matrix  $\Big(\mu_{\phi_{1,2}} (e_i,e_j) \Big)= s\mathbb I_{p,n}$ of $\mu_{\phi_{1,2}}$. In the basis  $(e_1,\dots,e_n, e_1^*,\dots,e_n^*)$ of $T^*\mathcal G$, where $(e_1^*,\dots,e_n^*)$ is the dual basis of $(e_1,\dots,e_n)$, the matrix $[\mu]$ of $\mu$ reads 
\begin{eqnarray}\label{eq:biinvariantmetricTsl(2,1)} [\mu]:=\begin{pmatrix}\mu(e_i,e_j) &\mu(e_i,e_j^*)\\ \mu(e_i,e_j^*) &\mu(e_i^*,e_j^*) \end{pmatrix}=
\begin{pmatrix}s\mathbb I_{p,q} & t\mathbb I_{n} \\ t\mathbb I_{n}  &{\mathbf 0}_n \end{pmatrix}\;.
\end{eqnarray}
The characteristic polynomial of $[\mu]$ is
$P(X) 
=(X+\lambda_1)^{p}(X+\lambda_2)^{p}(X-\lambda_1)^{n-p}(X-\lambda_2)^{n-p}$,
with $\lambda_1:=\frac{1}{2}(s-\sqrt{s^2+4t^2})$ and $\lambda_2:=\frac{1}{2}(s+\sqrt{s^2+4t^2})$.  So $[\mu]$  has $n$ positive  eigenvalues, $p$ of which being equal to $-\lambda_1$ and $n-p$ of which being equal to $\lambda_2$, and $n$ negative eigenvalues,
 $p$ of which coinciding with $-\lambda_2$ and $n-p$ of which being equal to $\lambda_1$. Hence, {\color{black} this} confirms that $\mu$ is of signature $(n,n)$. 

(B) If $n$ is  even, by Lemma \ref{ad-invariant-endomorphisms},  $\psi$, $\phi_{1,1}$, $\phi_{2,2}$ are of the forms $\psi=s_1\mathbb I_{\mathcal G} +s_2 J$,   $\phi_{1,1}=t_1\mathbb I_{\mathcal G} +t_2 J$  and  $\phi_{2,2}=t_1\mathbb I_{\mathcal G^*}+t_2J^T$,  for some fixed $J$ satisfying $J[x,y]=[Jx,y]$, for any $x,y\in\mathcal G $ and $s_1,s_2,t_1,t_2\in\mathbb R$ {\color{black} and $J^2=-\mathbb I_{\mathcal G}.$}  Thus, we have for any $x,y\in\mathcal G$, $f,g\in\mathcal G^*$, 
\beqn 
~~~~~~~\mu\Big((x,f),(y,g)\Big)&=&\mu(x,y)+\mu(x,g)+\mu(f,y)\nonumber\\
 &=& s_1K_0(x,y)+s_2K_0(J(x),y) +t_1\langle x,g\rangle\nonumber\\
&&+ t_2\langle J(x),g\rangle +t_1\langle y,f\rangle+ t_2\langle J(y),f\rangle 
\nonumber\\
&=&s_1K_0(x,y)+s_2K_J(x,y) +t_1\langle (x,f),(y,g)\rangle \cr 
& &+ t_2\langle (x,f),(y,g)\rangle_J . \hfill ~~~~~~~~~~~~~~~~~~~~~~~~~~~~~~~~\qed 
\nonumber
\eeqn

\section{Dual quaternions, dual split quaternions, rigid motions, screws motions}\label{chap:rigid-motions-dual} 
\subsection{A short tour on dual quaternions, dual split quaternions}\label{chap:reminders-quaternions-plit-quaternions}
\subsubsection{Quaternions, rigid motions, screw motions, dual quaternions}

Consider the following three vectors ${\mathbf i}$, ${\mathbf j}$, ${\mathbf k}$ called the pure imaginary units and satisfying the relation ${\mathbf i}^2={\mathbf j}^2={\mathbf k}^2={\mathbf i}{\mathbf j}{\mathbf k}=-1$ {\color{black}(along with associativity)}. 
A quaternion is of the form 
$Q=q_0+q_x{\mathbf i}+q_y{\mathbf j}+q_z{\mathbf k} =q_0+\vec{q}$, where $q_0,q_x,q_y,q_z\in\mathbb R$  and  $\vec{q}:=q_x{\mathbf i}+q_y{\mathbf j}+q_z{\mathbf k}$ is called the vector part or the pure imaginary part of $Q$. One $\mathbb R$-linearly extends the above relation to define a product on the space $\mathbb H$ of quaternions, given for any $P=p_0+\vec{p}$, $Q=q_0+\vec{q}\in\mathbb H$, by 
\beqn \label{Eq:productquaternions} 
PQ =q_0p_0 - \vec{q}\cdot\vec{p}+ q_0\vec{p} + p_0\vec{q} + \vec{q}\times\vec{p},
\eeqn
 where $ \vec{q}\cdot\vec{p}$ and  $ \vec{q}\times\vec{p}$ are respectively  the ordinary dot and vector products of $ \vec{q}$ and $\vec{p}$ seen as vectors in $\mathbb R^3$. A vector $u$ with coordinates $(u_x,u_y,u_z)$ in the canonical basis of $\mathbb R^3$, is identified with the pure  quaternion $\vec{u}=u_x{\mathbf i}+u_y{\mathbf j}+u_y{\mathbf k}$.
The product (\ref{Eq:productquaternions}) makes $\mathbb H$ an associative $\mathbb R$-algebra.
The conjugate of $Q$ is the quaternion
$Q^*=q_0-q_x{\mathbf i}-q_y{\mathbf j}-q_z{\mathbf k} =q_0-\vec{q}$.
As one can see $QQ^*=q_0^2+q_x^2+q_y^2+q_z^2=\vert\vert Q\vert\vert^2$, thus $Q$ is invertible if and only if $\vert\vert Q\vert\vert\neq 0$ and in this case, $Q^{-1} =\frac{1}{\vert\vert Q\vert\vert^2} Q^*$. The set of invertible quaternions is an analytic 4-dimenional Lie group, of whom the subset $\mathbb H_1$ of  those $Q$ with $\vert\vert Q\vert\vert=1$ (unit quaternions) is a closed subgroup, hence a Lie subgroup (according to Cartan's theorem). 
Consider  the $8$-dimensional $\mathbb R$-algebra $\mathfrak{gl}(2,\mathbb C)$ of $2\times 2$ matrices with coefficients in $\mathbb C$ and its $4$-dimensional $\mathbb R$-subalgebra  $\mathcal H$ spanned by the $2\times 2$ matrices 
$X_0:=E_{1,1}+E_{2,2}$, $X_1:= \sqrt{-1}\: (E_{2,2}-E_{1,1})$, $ X_2:=E_{2,1}-E_{1,2}$, $ X_3:=\sqrt{-1}\;(E_{1,2}+E_{2,1})$.
The algebra  $\mathbb H$ is isomorphic to  $\mathcal H$,
  via the linear map $\psi$, defined by 
\beqn \label{isomorphimSU(2)-quaternions}\psi({\mathbf 1})=X_0, \; \psi({\mathbf i})=X_1, \; \psi({\mathbf j})=X_2, \; \psi({\mathbf k})=X_3.
\eeqn   
 That is,  $\det \psi\neq 0$ and $ \psi(Q_1Q_2)=\psi(Q_1)\psi(Q_2)$, for any $Q_1,Q_2\in\mathbb H$. 
Given a unit quaternion $Q = q_0+ \vec{q}$, one sets $\theta = 2\cos^{-1}(q_0)$,  $0\leq \theta\le 2\pi$, and $ \vec{u} =\frac{\vec{q}}{\sin(\frac{\theta}{2})}$ if $\theta \neq 0,\;2\pi$ and  $ \vec{u} =0$ otherwise, so that 
$Q=\cos(\frac{\theta}{2})+\vec{u}\sin(\frac{\theta}{2})$, with $\vert\vert \vec{ u}\vert\vert=1$ {\color{black}when $\theta \neq 0,\;2\pi$.}
Note that $\mathcal H$ contains SU(2) and the restriction $\psi: \mathbb H_1\to SU(2)$  is an isomorphism between the Lie groups 
$\mathbb H_1$ and $SU(2)$.
More precisely,  let $Q=\cos(\frac{\theta}{2})+\sin(\frac{\theta}{2}) \vec{u} \in\mathbb H_1$, 
$\vec{u} = u_x{\mathbf i} +  u_y{\mathbf j}+ u_z{\mathbf k}$ with $u_x^2+u_y^2+u_z^2=1$, then  the image 
$\psi(Q)=\sigma$ is
\begin{eqnarray}
\psi(Q)= \begin{pmatrix}
\cos(\frac{\theta}{2})-iu_x\sin(\frac{\theta}{2}) & -u_y\sin(\frac{\theta}{2})+ iu_z\sin(\frac{\theta}{2}) \\
                                                                         &         \\
 u_y\sin(\frac{\theta}{2}) + iu_z\sin(\frac{\theta}{2}) &\cos(\frac{\theta}{2})+iu_x\sin(\frac{\theta}{2})
 \end{pmatrix} =
 \begin{pmatrix}
 z&-\bar w\\ 
   &    \\
 w&\bar z
 \end{pmatrix}, 
\nonumber
  \end{eqnarray} 
with $z=\cos(\frac{\theta}{2})-iu_x\sin(\frac{\theta}{2})$,  $w=u_y\sin(\frac{\theta}{2}) + iu_z\sin(\frac{\theta}{2})$,  $\vert z\vert^2 + \vert w\vert^2 =1$. {\color{black} Note that $\psi$ sends the vector part of $\mathbb H$ isomorphically to the Lie algebra $\mathfrak{su}(2).$}
On the other hand, one constructs a Lie group double covering, more precisely, the well known  Lie group homomorphism $\tilde \Pi: SU(2)\to SO(3)$ with kernel $\ker(\tilde \Pi)=\{1,-1\}$, using  
 the adjoint representation Ad: $SU(2)\to GL(\mathfrak{su}(2))$, of $SU(2)$ with Ad$_\sigma X=\sigma X\sigma^{-1} 
$. Indeed, the metric defined by $K(M,N):=$Trace$(MN)$  is negative definite and Ad-invariant on the Lie algebra $\mathfrak{su}(2)$ of $SU(2)$, as $K(Ad_\sigma X,Ad_\sigma Y)=K(X,Y)$. Thus,  $Ad_\sigma$ lies in $SO(3)$, for any $\sigma\in SU(2)$ and $\ker(Ad)=\{1,-1\}$. We set $Ad=:\tilde \Pi$.  
 Note that $(X_1,X_2,X_3)$ above is a basis of the Lie algebra $\mathfrak{su}(2)$ of  $SU(2)$.
For $Q=\cos(\frac{\theta}{2})+\vec{u}\sin(\frac{\theta}{2})\in\mathbb H_1$, the matrix of $\tilde \Pi\circ\psi (Q) $ in the basis  $(X_1,X_2,X_3)$, is
{\footnotesize
\begin{eqnarray}\label{MatrixAdSU(2)}
\begin{pmatrix} 
 \cos\theta+u_x^2(1-\cos\theta) &-\!u_z\sin\theta\!+\! u_xu_y(1\!-\!\cos\theta)  &(1\!-\!\cos\theta)u_xu_z\!+\! u_y\sin\theta\\
                                                                         &         \\
 u_x u_y (1\!-\!\cos\theta)\!+\!u_z \sin\theta &u_y^2(1\!-\!\cos\theta)\!+\! \cos\theta &u_yu_z(1\!-\!\cos \theta)\!-\!u_x\sin\theta \\
                                                &                            \\
u_x u_z (1\!-\!\cos\theta)\!-\!u_y \sin\theta  &  u_y u_z (1\!-\! \cos\theta)\!+\! u_x \sin\theta &  \cos\theta\!+\!u_z^2(1\!-\!\cos\theta)
\end{pmatrix}
\end{eqnarray} }
which is exactly the matrix of the rotation of angle $\theta$ around the axis defined by the unit vector $\vec{u}$, where the vector space underlying $\mathfrak{so}(3)$ has been identified with $\mathbb R^3$. 

 For a given unit $Q\in\mathbb H_1$, consider the linear map 
\beq \label{definition-RQ} 
\mathcal R_Q:\mathbb R^3\to\mathbb R^3,\;\; \mathcal R_Q(\vec{u})=Q\vec{u}Q^{-1}. 
\eeq
Each $\mathcal R_Q$ is a rotation of $\mathbb R^3$. 
More explicitly one has 
 $Q \vec{u} = ( - \vec{q}\cdot\vec{u}, q_0\vec{u} + \vec{q}\times\vec{u} )$ and
 \begin{eqnarray}\label{rotation-unit-quaternion1} Q \vec{u}Q^{-1} =
  (q_0^2+F(\vec{q}))^2\vec{u}
 + (\vec{q}\cdot\vec{u})\vec{q}  
= 
 (q_0^2 - \vert\vert\vec{q}\vert\vert^2)\vec{u} + 2q_0\vec{q}\times\vec{u} +2 (\vec{q}\cdot\vec{u})\vec{q} \;,
\end{eqnarray} 
where $F(x): \mathbb R^3\to \mathbb R^3$, $y\mapsto F(x)y= x\times y, $ is the cross product by $x$. So $F: \mathbb R^3\to \mathfrak{so}(3)$, $x\mapsto F(x)$, is an  isomorphism of vector spaces.
If the tensor product $x\otimes y$ is seen as a linear map whose matrix has coefficients $(x\otimes y)_{i,j} = x_iy_j$, then
  \begin{eqnarray} 
  (\vec{q}\cdot\vec{u})\vec{q}&=&q_x(q_xu_x+q_yu_y+q_zu_z) e_1+q_y(q_xu_x+q_yu_y+q_zu_z) e_2
\nonumber\\
&&+q_z(q_xu_x+q_yu_y+q_zu_z) e_3=(\vec{q}\otimes \vec{q}) \vec{u}.\end{eqnarray} 
So we get $Q \vec{u}Q^{-1} =\Big((q_0+F(\vec{q}))^2 + \vec{q}\otimes\vec{q} \Big)\vec{u}  $. So the matrix of the linear map $\mathcal R_Q=(q_0+F(\vec{q}))^2 + \vec{q}\otimes\vec{q} $, in the canonical
basis of $\mathbb R^3$, is

$\begin{pmatrix}
q_0^2+q_x^2-q_y^2-q_z^2 & -2q_0q_z+2q_xq_y & 2q_0q_y+2q_xq_z\\ 
                        &          &   \\
2q_0q_z+2q_xq_y & q_0^2-q_x^2+q_y^2-q_z^2 & -2q_0q_x+2q_yq_z\\ 
   &          &   \\
-2q_0q_y+2q_xq_z & 2q_0q_x+2q_yq_z & q_0^2-q_x^2-q_y^2+q_z^2 
\end{pmatrix}$,

\noindent hence $\mathcal R_Q\mathcal R_Q^T=\vert\vert Q\vert\vert^2\; \mathbb I_{\mathbb R^3}=\mathbb I_{\mathbb R^3}$, which confirms that $\mathcal R_Q$ is a rotation in $\mathbb R^3$.
  If one sets $Q=\cos(\frac{\theta}{2})+\vec{\omega} \sin(\frac{\theta}{2})$, where  $\vec{\omega}=\omega_x{\mathbf i}+\omega_y{\mathbf j}+\omega_y{\mathbf k}$,  $\omega_x^2+\omega_y^2+\omega_z^2=1$ and $0\leq\theta\leq 2\pi$, then the matrix of $\mathcal R_Q$ now reads
  {\footnotesize
 \begin{eqnarray}\label{Ad2}
\begin{pmatrix}\cos\theta+\omega_x^2(1-\cos\theta) &\omega_x\omega_y(1-\cos\theta)-\omega_z\sin\theta &\omega_x\omega_z(1-\cos\theta)+\omega_y\sin\theta\\
   &          &   \\
\omega_x\omega_y(1-\cos\theta)+\omega_z\sin\theta & \cos\theta+\omega_y^2(1-\cos\theta) &\omega_y\omega_z(1-\cos\theta)-\omega_x\sin\theta\\
   &          &   \\
\omega_x\omega_z(1-\cos\theta)-\omega_y\sin\theta &\omega_y\omega_z(1-\cos\theta)+\omega_x\sin\theta&\cos\theta+\omega_z^2(1-\cos\theta)\end{pmatrix}.
\end{eqnarray} }
Clearly, this corresponds to the matrix of the rotation with angle $\theta$ around the axis defined by the unit vector $\omega$ with components $(\omega_x,\omega_y,\omega_z)$ in the canonical basis of $\mathbb R^3$.
So $\mathcal R:\mathbb H_1\to SO(3), Q\mapsto \mathcal R_Q$, is a Lie group homomorphism, with kernel $\{1,-1\}$. Indeed one has 
$\mathcal R_{Q_1Q_2}(\vec{u})=Q_1Q_2uQ_2^{-1} Q_1^{-1}=\mathcal R_{Q_1}\mathcal  R_{Q_2}(\vec{u})$ for any $Q_1,Q_2\in \mathbb H_1$ and any $u\in \mathbb R^3$. That is, $\mathcal R_{Q_1Q_2} = \mathcal R_{Q_1}\mathcal R_{Q_2}$.  On the other hand,  suppose that $\mathcal R_{Q} =\mathbb I_{\mathbb R^3}$. Then choosing $\vec{u}$ perpendicular to $\vec{q}$ 
so that  $\vec{u}$ and $\vec{q}\times\vec{u}$ are linearly independent and  $\vec{q}\cdot\vec{u}=0$,  implies that the expression 
$\vec{u}=Q \vec{u}Q^{-1} = (q_0^2 - \vert\vert\vec{q}\vert\vert^2)\vec{u} + 2q_0\vec{q}\times\vec{u} +2 (\vec{q}\cdot\vec{u})\vec{q} $ becomes 
$\vec{u}= (q_0^2 - \vert\vert\vec{q}\vert\vert^2)\vec{u} + 2q_0\vec{q}\times\vec{u} $, or equivalently, $\vec{q}=0$ and $q_0^2=1$. This is,  $Q=\pm 1\;. $  
Conversely,  to any given rotation $\hat{\mathcal R}$ in $3$D, of angle $\theta$  around an axis of rotation defined by a chosen unit vector, say  $\omega\in\mathbb R^3$, one associates a unit quaternion 
$Q =\cos(\frac{\theta}{2}) +\omega\sin(\frac{\theta}{2})$. Further applying  to the corresponding rotation $\mathcal R_Q$, the expression $Q \vec{u}Q^{-1} = (q_0^2 - \vert\vert\vec{q}\vert\vert^2)\vec{u} + 2q_0\vec{q}\times\vec{u} +2 (\vec{q}\cdot\vec{u})\vec{q} $ obtained in (\ref{rotation-unit-quaternion1}),  leads to the following equalities, valid for any $\vec{u}:$ 
 \begin{eqnarray}\mathcal R_Q(\vec{u}) &=&
\Big( \cos\theta\; + \sin\theta\; F(\vec{\omega}) +(1-\cos\theta)\;\vec{\omega}\otimes \vec{\omega} \Big)\;\vec{u}
=\hat{\mathcal R}(\vec{u}) \;.\end{eqnarray}
 So $\hat{\mathcal R}=\mathcal R_Q$, and hence $\mathcal R$ is surjective. Consider the coset  $\mathbb H_1\slash \{+1,-1\}$ of elements $\bar Q=\{Q,-Q\}=Q\mod \pm1$.  From the above, the map 
$\bar {\mathcal R}:\mathbb H_1\slash \{+1,-1\}\to SO(3)$, $\bar Q\mapsto \mathcal R_Q$, is a Lie group isomorphism. 
Thus $\mathbb H_1$ is the universal cover of SO($3$)  and  $\mathcal  R$  
 can be considered as the corresponding universal covering. We note that, comparing (\ref{MatrixAdSU(2)}) and (\ref{Ad2}), one sees that
 \begin{eqnarray}\label{equality}\tilde \Pi\circ \psi = \mathcal R \;. \end{eqnarray}
One naturally extends $\tilde \Pi$ to a surjective homomorphism  $T^*SU(2)\to T^*SO(3),$ with kernel  $\{+1,-1\}.$

Dual quaternions are defined as combinations of ordinary quaternions  and dual numbers. More precisely, a dual quaternion is of the form 
 $\hat Q={Q}_r+ \mathfrak{e}{Q}_d$ where ${Q}_r$ and ${Q}_d$ are quaternions and $\mathfrak{e}$, termed the dual operator, satisfies the properties
$\mathfrak{e}\neq 0$ and $\mathfrak{e}^2=0$. 
We have the following rules.
The sum  and product of two dual quaternions $\hat Q_1=Q_{r,1} + \mathfrak{e}Q_{d,1} $ and $\hat Q_2=Q_{r,2} + \mathfrak{e}Q_{d,2}$, are respectively 
\begin{eqnarray}\label{sum-product-dual-quaternions}
\hat Q_1+\hat Q_2&=&Q_{r,1}+Q_{r,2} +  \mathfrak{e}(Q_{d,1}+Q_{d,2}) \; ,\nonumber\\ 
\hat Q_1\hat Q_2&=&Q_{r,1}Q_{r,2} + \mathfrak{e}(Q_{r,1}Q_{d,2}+Q_{d,1}Q_{r,2}).
\end{eqnarray}
Dual quaternions were introduced by Clifford in 1882 (see \cite{clifford}). They have since gained  a lot of interest due to their applications in many areas such as Robotics, Computer graphics and video games. 
They appear as a natural unified framework for both the translations and rotations, with  computational advantages  such as being robust and presenting singularity free solutions.

 The dual conjugation (analogous to complex conjugation) of 
$\hat{Q}$ is 
$\hat{Q}^*={Q}_r^*+  \mathfrak{e}{Q}_d^*$  and its dual quaternionic conjugation (analogous to quaternionic conjugation) is
$\bar{\hat{Q}}={Q}_r-  \mathfrak{e}Q_d$. 
One also defines the magnitude of a dual quaternion $\hat{Q}={Q}_r+ \mathfrak{e}{Q}_d$ using the expression  
$\vert\vert \hat{Q} \vert\vert^2 :=\hat{Q}\hat{Q}^*= {Q}_{r}{Q}_{r}^* + \mathfrak{e}({Q}_{r}{Q}_{d}^*+{Q}_{d}{Q}_{r}^*) $.
If we set ${Q}_r=s_r+x_r{\mathbf i} +y_r{\mathbf j}+z_r{\mathbf k}$, ${Q}_d=s_d+x_d{\mathbf i} +y_d{\mathbf j}+z_d{\mathbf k}$, we can see the following
\begin{eqnarray} {Q}_r^*{Q}_d+{Q}_d^*{Q}_r&=&
2(s_rs_d+x_rx_d+y_ry_d+z_rz_d)=2\langle{Q}_r,{Q}_d\rangle\;.\end{eqnarray}
So we have $\vert\vert \hat{Q} \vert\vert^2 =\vert\vert {Q}_r \vert\vert^2 +2 \mathfrak{e}\langle{Q}_r,{Q}_d\rangle $.
 We also further have
\begin{eqnarray} \label{norm-dual-quaternion-multiplicative}\vert\vert\hat Q_1\hat{Q}_2\vert\vert^2&=&
\vert\vert \hat{Q}_1 \vert\vert^2 \; \vert\vert \hat{Q}_2 \vert\vert^2.
 \end{eqnarray}

A dual quaternion $\hat Q:=Q_r+\mathfrak{e}Q_d$ has an inverse if and only if $Q_r\neq 0$. In that case, $\hat Q^{-1}= Q_r^{-1}-\mathfrak{e} (Q_r^{-1}Q_dQ_r^{-1})$.
With the sum and the  product (\ref{sum-product-dual-quaternions}), the $8$-dimensional vector space underlying the space $\hat{\mathbb H}$ of dual quaternions is an associative algebra with a unit, namely $1$. Hence its (open) subset consisting of invertible dual quaternions is an analytic Lie group of dimension $8$, (see e.g. J. P. Serre \cite{serre}, p. 103, for the general case of the set of invertible elements of an associative algebra with unit)  of which the subset $\hat{\mathbb H}_1$ of unit dual quaternions is a closed subgroup and hence a closed Lie subgroup, according to Cartan's Theorem. 
Another useful way of representing unit dual quaternions is as follows 
$\hat Q=Q+\mathfrak{e}(\vec{u} Q)$, where $Q$ is a unit quaternion,  $\vec{u}$ a pure quaternion. 
From e.g. \cite{sastry}, 'a rigid motion of an object (in the Euclidean 3-space $\mathbb R^3$) is a continuous movement of the particles
in the object such that the distance between any two particles remains
fixed at all times. A rigid motion is one that preserves the distance
between points and the angle between vectors. The net movement of a rigid body from one location
to another via a rigid motion is called a rigid displacement. In general,
a rigid displacement may consist of both translation and rotation of the
object.'
 A rigid displacement is represented by an element of the special affine group SE($3$):=SO($3$)$\ltimes \mathbb R^3$, that is, by an affine map $x\mapsto L(x)+v$, with linear part $L\in SO(3)$ and fixed vector $v\in\mathbb R^3$.   A rigid motion is represented by a curve on SE($3$).
It is well known that, together with the composition of maps, SE($3$) is a Lie group. A screw displacement is a rigid displacement consisting of a rotation at constant angular velocity around an axis (the screw or twist axis) followed by a translation with constant  (translationnal) velocity along the same axis. 
{\color{black}From Chasles theorem, any rigid motion can be realized as a screw motion.} 
One extends the map $\mathcal R$ discussed above, to {\color{black} the} map $\Pi: \hat{\mathbb H}_1\to SE(3) $ as follows,
\begin{eqnarray}\label{formula:unit-dual-quaternions2a}\Pi: \hat{\mathbb H}_1\to SE(3) \;,\; \hat Q=Q+\mathfrak{e}(\vec{u} Q)  \;\mapsto \; (\mathcal  R_Q, \vec{u}). \end{eqnarray}

If $\hat Q=Q+\mathfrak{e}\vec{u} Q $ and $\hat P=P+\mathfrak{e}\vec{v} P $  are unit dual quaternions, then 
the unit dual quaternion $\hat Q\hat P$ is represented as $\hat Q\hat P= QP+ \mathfrak{e}(Q\vec{v}P +\vec{u}QP ) = QP+ \mathfrak{e}(Q\vec{v}Q^{-1} +\vec{u})QP   = QP+ \mathfrak{e}(\vec{X}QP) $ with $\vec{X}   = \vec{u}+\mathcal  R_Q\vec{v}$.
 Thus, we have
$\Pi(\hat Q\hat P)=  (\mathcal R_{QP}, \vec{u}+\mathcal  R_Q\vec{v})$, which we rewrite as 
 $\Pi(\hat Q\hat P)= (\mathcal R_Q, \vec{u})(\mathcal R_P, \vec{v})= \Pi(\hat Q)\Pi(\hat P)$.
A unit dual quaternion $\hat Q:=Q+\mathfrak{e}(\vec{u} Q)$ satisfies $\Pi(\hat Q) =(\mathbb I_{\mathbb R^3},0)$ if and only if 
$\mathcal  R_Q=\mathbb I_{\mathbb R^3}$ and $\vec{u}=0$,  or equivalently $\hat Q=\pm 1$. Of course, for any rigid displacement $(\hat{\mathcal  R},v)$, where $\hat {\mathcal R}$ is a rotation of angle $\theta$ around an axis defined by a unit vector $u$, the unit dual  quaternion  $\hat Q=Q+\mathfrak{e}\vec{v} Q $, where $Q=\cos(\frac{\theta}{2})+\sin(\frac{\theta}{2})\;\vec{u}$, satisfies $\Pi(\hat Q) =(\hat {\mathcal R},v)$.  So $\Pi$ is a surjective homomorphism between the Lie groups $\hat{\mathbb H}_1$ and $SE(3)$, with  $\ker\Pi=\{1,-1\}$, hence a double cover. As $\hat{\mathbb H}_1$ is simply connected,  $\Pi: \hat{\mathbb H}_1\to$ SE($3$) can be seen as the universal covering of SE($3$).
This allows one to represent rigid motions by dual quaternions.

Now our Theorem \ref{theorem:commutingdiagramSO(3)} asserts that, up to isomorphism, the Lie group $\hat{\mathbb H}_1$ is in fact the cotangent bundle $T^*SU(2)$ of the compact simple Lie group SU(2), endowed with its Lie group structure induced by right trivialization. Theorem \ref{theorem:commutingdiagramSO(3)} further shows that the group SE(3) of Rigid motions of the Euclidean space is the cotangent bundle of  the compact simple Lie group $SO(3)$, up to isomorphism of Lie groups.

\subsubsection{Split quaternions, rigid motion in Minkowski space, dual split quaternions}
A split quaternion is of the form 
$Q=q_0+q_x{\mathbf i}+q_y{\mathbf j}+q_z{\mathbf k} =q_0+\vec{q}$, where $q_0,q_x,q_y,q_z\in\mathbb R$  and  
$\vec{q}:=q_x{\mathbf i}+q_y{\mathbf j}+q_z{\mathbf k}$ is called the vector part or the pure imaginary part of $Q$.  The vectors ${\mathbf i}$, ${\mathbf j}$, ${\mathbf k}$  satisfy ${\mathbf i}^2=-1$, ${\mathbf j}^2={\mathbf k}^2=1$ and ${\mathbf i}{\mathbf j}{\mathbf k}=1$ {\color{black}(along with associativity)}. One $\mathbb R$-linearly extends the latter relations to define a product in the space of split quaternions, hereafter denoted by $\mathbb S$, as follows.
 If $P=p_0+p_x{\mathbf i}+p_y{\mathbf j}+p_z{\mathbf k}=p_0+\vec{p} $ is another split quaternion, the product $Q P$ of $Q$ and $P$, is 
{\color{black}\begin{eqnarray} \label{eq:product-split-quaternions1}QP
&=& q_0p_0 - (p_xq_x - p_yq_y - p_zq_z) +q_0\vec{p}+p_0\vec{q}\nonumber\\
&&+(q_zp_y-q_yp_z){\mathbf i}
+(q_zp_x-q_xp_z){\mathbf j}+(q_xp_y-q_yp_x){\mathbf k}.
\end{eqnarray}
}
One identifies $\mathbb S$ with the vector space  $\mathbb R^4$, where a split quaternion $Q=q_0+q_x{\mathbf i}+q_y{\mathbf j}+q_z{\mathbf k} $ is identified with the vector with coordinates $(q_0,q_x,q_y,q_z)$ in the canonical basis of $\mathbb R^4$. 
One also considers  the metric $\langle , \rangle$ on $\mathbb R^4$ with signature $(2,2)$ defined by
$\langle Q, P\rangle =q_0p_0+q_xp_x - p_yq_y - p_zq_z$,  the restriction to $\mathbb R^3$  of which is
 the Lorentzian metric given by  $\langle \vec{q}, \vec{p}\rangle:=p_xq_x - p_yq_y - p_zq_z$. 
We will  let $\mathbb R_1^3$ stand for  $\mathbb R^3$ together with the Lorentzian metric $\langle,\rangle$, whereby the canonical basis of $\mathbb R^3$ is identified with $({\mathbf i, \; \mathbf j, \; \mathbf k})$, and refer to it as the   Minkowski  $3$-space.
 Let us consider the cross product $\times_s$ in the Minkowski  $3$-space, defined for $\vec{q}=q_x{\mathbf i}+q_y{\mathbf j}+q_z{\mathbf k} $
 and $\vec{p}=p_x{\mathbf i}+p_y{\mathbf j}+p_z{\mathbf k} $ as,
\begin{eqnarray} \label{cros-product-Minkowski}\vec{q}\times_s\vec{p}:=(q_zp_y-q_yp_z){\mathbf i}
+(q_zp_x-q_xp_z){\mathbf j}+(q_xp_y-q_yp_x){\mathbf k}.
\end{eqnarray}
The product (\ref{eq:product-split-quaternions1}) can then be rewritten as 
{\color{black}\begin{eqnarray} \label{eq:product-split-quaternions2}QP &=&q_0p_0-\langle \vec{q}, \vec{p}\rangle
+q_0\vec{p}+p_0\vec{q}+\vec{q}\times_s\vec{p}.
\end{eqnarray}
}
From (\ref{eq:product-split-quaternions2}), one notes that  $QP-PQ=\vec{q} \vec{p} - \vec{p} \vec{q} = 2\vec{q}\times_s\vec{p}$, so one can also define the cross product $\times_s$ as 
\begin{eqnarray} \label{eq:cross-product-split-quaternions3}\vec{q}\times_s\vec{p}:=\frac{1}{2}(\vec{q} \vec{p} - \vec{p} \vec{q}),
\end{eqnarray}
where $\vec{q} \vec{p} =- p_xq_x + p_yq_y + p_zq_z +(q_zp_y-q_yp_z){\mathbf i}
+(q_zp_x-q_xp_z){\mathbf j}+(q_xp_y-q_yp_x){\mathbf k}$ is the product of $\vec{q} $ and $ \vec{p} $.
Like in the Euclidean case, one also has $\langle\vec{q},\vec{q}\times_s \vec{p}\rangle = \langle\vec{p},\vec{q}\times_s \vec{p}\rangle=0$ for any $\vec{q}, \vec{p}\in\mathbb R^3_1$.
The conjugate of a split quaternion $Q$ is defined as 
$Q^*=q_0-q_x{\mathbf i}-q_y{\mathbf j}-q_z{\mathbf k} =q_0-\vec{q}$. Considering the product $QQ^*= q_0^2+q_x^2 -q_y^2 - q_z^2 =\langle Q, Q\rangle$, one sees that  a split quaternion $Q$ admits an inverse if and only if $\langle Q, Q\rangle\neq 0$ and in that case $Q^{-1} =\frac{1}{\langle Q, Q\rangle} Q^*$.
We say that $Q$ is a unit split quaternion, if $\langle Q,Q\rangle =1$. We let $\mathbb S_1$ stand for the (closed) subgroup of $\mathbb S$ made of unit split quaternions. A pure split quaternion $\vec{q}$ is said to be timelike if $\langle\vec{q},\vec{q}\rangle> 0$,    lightlike if $\langle\vec{q},\vec{q}\rangle =0$ and  $\vec{q}\neq 0$,  
and spacelike if  $\langle\vec{q},\vec{q}\rangle < 0$  or $\vec{q}=0$. 
 From (\ref{eq:product-split-quaternions2}), we see that any pure split quaternion  $\vec{q}$ satisfies $\vec{q} ^2= - \langle \vec{q}, \vec{q}\rangle$.
 Let $Q=q_0+\vec{q}$ be a unit split quaternion, then either $Q=\pm 1$ or $\langle \vec{q}, \vec{q}\rangle \neq  0$. In the case $\langle \vec{q}, \vec{q}\rangle \neq  0$, one sets $\vec{u}=\frac{\vec{q}}{\sqrt{\vert \langle\vec{q},\vec{q}\rangle\vert}}$. 
 So if  $\langle\vec{q},\vec{q}\rangle\geq 0$, one sets  $\theta=\cos^{-1} q_0$, with $0\leq \theta\leq \pi$, so that $q_0=\cos(\theta)$ and $\sin(\theta)=\sqrt{\langle\vec{q},\vec{q}\rangle}$. Furthermore, in the case  $\langle\vec{q},\vec{q}\rangle> 0$,   the unit  pure split quaternion   
$\vec{u}=\frac{\vec{q}}{\sqrt{\langle\vec{q},\vec{q}\rangle}}$  satisfies $\vec{u}^2= -1$ and $Q=q_0+\vec{q}=\cos\theta+\sin\theta\;\vec{u}
=\displaystyle\sum_{k=0}^\infty\frac{(\theta\vec{u})^{2k}}{(2k)!}+\displaystyle\sum_{k=0}^\infty\frac{(\theta\vec{u})^{2k+1}}{(2k+1)!}=\displaystyle\sum_{k=0}^\infty\frac{(\theta\vec{u})^k}{k!}=e^{\theta\vec{u}}$.                        
So a unit split quaternion $Q$ with a timelike vector part can always be written in the form $Q=\cos\theta+\sin\theta\;\vec{u}=e^{\theta\vec{u}}$, with $0\leq \theta\leq\pi$ and $\vec{u}$ is a unit (timelike) pure split quaternion. 
Now if  $\langle\vec{q},\vec{q}\rangle < 0$, one sets $\gamma=\cosh^{-1} q_0$, the pure split quaternion   
$\vec{u}=\frac{\vec{q}}{\sqrt{-\langle\vec{q},\vec{q}\rangle}}$  satisfies $\vec{u}^2= 1=-\langle \vec{u}, \vec{u}\rangle$ and $Q=q_0+\vec{q}=\cosh\theta+\sinh\theta\;\vec{u}
=\displaystyle\sum_{k=0}^\infty\frac{(\theta\vec{u})^{2k}}{(2k)!}+\displaystyle\sum_{k=0}^\infty\frac{(\theta\vec{u})^{2k+1}}{(2k+1)!}=\displaystyle\sum_{k=0}^\infty\frac{(\theta\vec{u})^k}{k!}=e^{\theta\vec{u}}$.
Thus,  a unit split quaternion $Q$ with a spacelike vector part can be written in the form $Q=\cosh\theta+\sinh\theta\;\vec{u}=e^{\theta\vec{u}}$, with $0\leq \theta\leq\pi$ and $\vec{u}$ is a unit (spacelike) pure split quaternion. 
This allows to write any unit split quaternion as the exponential of some pure 
 split quaternion. For more details, see \cite{kula-yayli}, \cite{kula-yayli-screw-motion}, where rotations on the Minkowski space are also discussed using split quaternions.
Let $\mathcal M(2,\mathbb R)$ be the space of linear transformations of $\mathbb R^2,$ identified with that of $2\times 2$ real matrices (on the canonical basis). The well  known linear map 
\beqn \label{splitquaternions-as-2by2-matrices}
 \mathfrak{\omega}:\mathbb S\to  \mathcal M(2,\mathbb R),\; \mathfrak{\omega}(1)= E_{1,1}+ E_{2,2},  \;
\mathfrak{\omega}(\mathbf i)= E_{1,2}- E_{2,1},
\nonumber
\\
 \mathfrak{\omega}(\mathbf j)= E_{1,2}+ E_{2,1},\; 
 \mathfrak{\omega}(\mathbf k)= E_{1,1}- E_{2,2},
\eeqn
 is an isomorphism of associative algebras with units,  and 
$\langle  Q, Q\rangle=\det\mathfrak{\omega}( Q).$ 
In particular, $\mathfrak{\omega}$ restricts to an isomorphism $\mathfrak{\omega}: \mathbb S_1\to SL(2,\mathbb R)$ between the Lie groups $\mathbb S_1$ and $SL(2,\mathbb R).$
{\color{black} Note that $\mathfrak{\omega}$ maps the vector part of $\mathbb S$ isomorphically to the Lie algebra $\mathfrak{sl}(2,\mathbb R).$}
Like for quaternions, {\color{black} we define} a Lie group homomorhism, also denoted by  $\mathcal R,$
\beqn \label{def:R} \mathcal R: \mathbb S_1\to SO(2,1)\;, \text{with } \mathcal R_Q:\mathbb R^3\to \mathbb R^3, \; \mathcal R_Q(\vec{v})=Q\vec{v} Q^{-1}.\eeqn
Similarly to dual quaternions, dual split quaternions are combinations of split quaternions  and dual numbers. A dual split quaternion is of the form 
 $\tilde Q={Q}_r+\mathfrak{e}{Q}_d$ where ${Q}_r$ and ${Q}_d$ are split quaternions and the dual operator $\mathfrak{e}$ satisfies
$\mathfrak{e}\neq 0$ and $\mathfrak{e}^2=0$, as above.
The following rules apply.
The sum  and product of two split dual quaternions $\tilde Q_1=Q_{r,1} + \mathfrak{e} Q_{d,1}$ and $\tilde Q_2=Q_{r,2} + \mathfrak{e}Q_{d,2}$, are respectively  defined as 
\begin{eqnarray}
\tilde Q_1+\tilde Q_2&=&Q_{r,1}+Q_{r,2} + \mathfrak{e} (Q_{d,1}+Q_{d,2}), \nonumber\\ 
\tilde  Q_1\tilde  Q_2&=&Q_{r,1}Q_{r,2} +  \mathfrak{e} (Q_{r,1}Q_{d,2}+Q_{d,1}Q_{r,2}). \label{sum-product-dual-split-quaternions}
\end{eqnarray}
Here too, the dual conjugation of 
$\tilde{Q}$ is 
$\tilde{Q}^*\!=\!{Q}_r^*\!+\!  \mathfrak{e} {Q}_d^*$  and its dual split quaternionic conjugaison  is
$\bar{\tilde{Q}}\!=\!{Q}_r\!-\!  \mathfrak{e}{Q}_d$. 
The magnitude of a dual split quaternion is defined using the expression  $\vert\vert \tilde{Q} \vert\vert^2 \!:=\!\tilde{Q}\tilde{Q}^*\!=\! {Q}_{r}{Q}_{r}^* \!+\!  \mathfrak{e} ({Q}_{r}{Q}_{d}^*\!+\!{Q}_{d}{Q}_{r}^*)$.
If we set ${Q}_r=s_r+x_r{\mathbf i} +y_r{\mathbf j}+z_r{\mathbf k}$, ${Q}_d=s_d+x_d{\mathbf i} +y_d{\mathbf j}+z_d{\mathbf k}$, we can see the following
\begin{eqnarray} 
{Q}_r^*{Q}_d+{Q}_d^*{Q}_r 
&=&2(s_rs_d+x_rx_d+y_ry_d+z_rz_d)=2\langle{Q}_r,{Q}_d\rangle.
\end{eqnarray}
So we also have $\vert\vert \tilde{Q} \vert\vert^2 =\vert\vert {Q}_r \vert\vert^2 +2  \mathfrak{e} \langle{Q}_r,{Q}_d\rangle$.
 Exactly in the same way as in (\ref{norm-dual-quaternion-multiplicative}), we also have
\begin{eqnarray}\label{norm-dual-split-quaternion-multiplicative}
\vert\vert\tilde Q_1\tilde{Q}_2\vert\vert^2&=&({Q}_{r,1}{Q}_{r,2}) ({Q}_{r,1}{Q}_{r,2} )^*+
2 \mathfrak{e} \langle{Q}_{r,1}{Q}_{r,2} \;,\;{Q}_{r,1}{Q}_{d,2}+{Q}_{d,1}{Q}_{r,2}\rangle \nonumber\\
&=&\vert\vert \tilde{Q}_1 \vert\vert^2 \; \vert\vert \tilde{Q}_2 \vert\vert^2.
 \end{eqnarray}
A dual split quaternion $\tilde Q:=Q_r+\mathfrak{e}Q_d$ has an inverse if and only if $\vert\vert Q_r\vert\vert \neq 0$. In that case, $\tilde Q^{-1}= Q_r^{-1}-\mathfrak{e} (Q_r^{-1}Q_dQ_r^{-1})$.
With the sum and the  product (\ref{sum-product-dual-split-quaternions}), the space $\tilde{\mathbb S}$ is an associative algebra with a unit, namely $1$. Hence its subset consisting of invertible dual split quaternions is an analytic Lie group and the subset $\tilde{\mathbb S}_1$ of unit dual split quaternions is a closed Lie subgroup of the latter. Like for dual quaternions, a unit dual split quaternion  $\tilde Q$, can be represented as follows 
$\tilde Q=Q+\mathfrak{e}(\vec{u} Q)$, where $Q$ is a unit split  quaternion,  $\vec{u}$ a pure split quaternion. 

Now consider the Lie group $SO(2,1)\subset$GL($3,\mathbb R$), of  oriented linear map which preserve the Lorentz metric in $\mathbb R_1^3.$
The Lie algebra $\mathfrak{so}(2,1)$ of $SO(2,1)$  can be represented as the following Lie algebra of $3\times 3$ matrices
\beq
 \Big\{n_3(E_{1,2}+E_{2,1})-n_2(E_{1,3}+E_{3,1})+n_1(E_{3,2}-E_{2,3}), \;
 n_1,n_2,n_3\in\mathbb R\Big\}.\nonumber\eeq
  The linear map $\mathfrak{f}: \mathfrak{so}(2,1)\to \mathfrak{sl}(2,\mathbb R),$  given on the basis $E_1:=E_{1,2}+E_{2,1},$  $E_2:=-(E_{1,3}+E_{3,1}),$ $E_3:=E_{3,2}-E_{2,3}$ of  $\mathfrak{so}(2,1),$
by 
\beqn\label{isomorphism-so(2,1)-sl(2)}
\mathfrak{f}(E_1)=\begin{pmatrix}\frac{\sqrt{5}}{4}& -\frac{1}{2}\\ \frac{1}{8}& -\frac{\sqrt{5}}{4}\end{pmatrix},
\; 
\mathfrak{f}(E_2)=-\begin{pmatrix}0& 1\\ \frac{1}{4}& 0\end{pmatrix},\;
\mathfrak{f}(E_3)=\begin{pmatrix}\frac{1}{4}& -\frac{\sqrt{5}}{2}\\ \frac{\sqrt{5}}{8}& -\frac{1}{4}\end{pmatrix},
\eeqn
 is an isomorphism between the Lie algebras $ \mathfrak{so}(2,1)$ and $ \mathfrak{sl}(2,\mathbb R).$

 Let us also consider the linear map
\begin{equation}\label{defineH}
H:\mathbb R^3\to\mathfrak{so}(2,1),  \: x\mapsto H(x),
\end{equation}
 where $H(x)$ is the cross-product by $x$ as in (\ref{eq:cross-product-split-quaternions3}), that is, 
\begin{equation}\label{defineH(x)} 
H(x): \mathbb R^3\to\mathbb R^3\;,\;\; y\mapsto H(x)y:=x\times_s y.
\end{equation}
 If $(n_1,n_2,n_3)$ are the coordinates of  ${\mathbf n}\in\mathbb R^3$ in the canonical basis of $\mathbb R^3$, then 
 $n_3(E_{1,2}+E_{2,1})-n_2(E_{1,3}+E_{3,1})+n_1(E_{3,2}-E_{2,3})$ is the matrix of $H({\mathbf n})$ in the same basis. 
Clearly, $H$ is an isomorphism between $\mathbb R^3$ and the vector space underlying $\mathfrak{so}(2,1)$.
Recall that the universal cover of $SL(2,\mathbb R)$ and SO(2,1), is not a Lie group of matrices. This is a well-known common fact for all special linear groups SL($n,\mathbb R$) for any $n\ge 2$. 
We consider the Lie subgroup $SE(2,1):=SO(2,1)\ltimes \mathbb R^3 $ of the Lie group Aff($3,\mathbb R$) $=GL(3,\mathbb R)\ltimes \mathbb R^3$, made  of invertible affine displacements of $\mathbb R^3$ whose linear parts are elements of SO(2,1). A natural way to represent the Lie group $SE(2,1)$  is as the following group of $4\times 4$ real matrices
\begin{eqnarray}\label{defineSE(2,1)} SE(2,1)=\Big\{\begin{pmatrix}A&v\\ 0&1\end{pmatrix}, \;\; A\in SO(2,1),\; v\in\mathbb R^3\Big\}.\end{eqnarray} 
By analogy with the group $SE(3):=SO(3)\ltimes \mathbb R^3$ of rigid displacement of the Euclidean $3$-space, we call  $SE(2,1)$ the group of rigid displacements of the Minkowski $3$-space. By a rigid motion in the Minkowski $3$-space, one means a curve in $SE(2,1)$. One also establishes {\color{black} the following surjective}  homomorphism {\color{black}(double cover) }
\begin{eqnarray}\label{isomorphism-unitquaternions-SE(2,1)}\tilde{\mathbb S}_1\to SE(2,1)\;,\;\;\tilde Q:=Q+ \mathfrak{e}\vec{u}Q\mapsto ({\mathcal  R}_Q,\vec{u}),\end{eqnarray} 
 between the Lie group $\tilde{\mathbb S}_1$ of unit dual split quaternions  and the Lie group SE($2,1$).
In a similar way as for dual quaternions, Theorem \ref{theorem:cummutingdiagramSL(2)} says that,  up to isomorphism, the Lie group $\tilde{\mathbb S}_1$
 is the cotangent bundle $T^* SL(2,\mathbb R)$ of the simple Lie group  $SL(2,\mathbb R)$, endowed with its Lie group structure induced by right trivialization. 

\subsection{Main results on dual quaternions, dual split quaternions}\label{chap:rigid-motions-dual} 

\begin{theorem}\label{theorem:commutingdiagramSO(3)}The Lie group $SE(3)$  of rigid motions of the Euclidean space $\mathbb R^3$, is isomorphic to both the cotangent bundle $T^*SO(3)$ and  tangent bundle $TSO(3)$ of $SO(3)$, endowed with their Lie group structure induced by the right trivialization.  The Lie group $\hat{\mathbb H}_1$ of unit dual quaternions  
 is isomorphic to both the cotangent bundle $T^*SU(2)$ and  tangent bundle $TSU(2)$ of $SU(2)$, endowed with their Lie group structure induced by the right trivialization. 
\end{theorem}
\begin{proof} 
Consider the linear map $F: \mathbb R^3\to \mathfrak{so}(3)$, $x\mapsto F(x)$, where $F(x)$ is the cross product  by $x$, that is, $F(x)y:= x\times y$. So $F$ is an isomorphism between $\mathbb R^3$ and the vector space underlying $\mathfrak{so}(3)$. 
Let $S:SO(3)\to SO(3)$, $\sigma\mapsto S_\sigma$, where $S(\sigma) x=F^{-1}(\sigma F(x)\sigma^{-1})$, $x\in\mathbb R^3$. Note that $S$ is a smooth map, for it only involves smooth maps, the linear map $F$ and the adjoint map. 
We further have 
$S(\sigma\tau)x=F^{-1}(\sigma\tau F(x)\tau^{-1}\sigma^{-1})=F^{-1}(\sigma F(F^{-1}(\tau F(x)\tau^{-1}))\sigma^{-1})
= F^{-1}(\sigma F(S(\tau)x)\sigma^{-1})= S(\sigma) S(\tau)x$,
and $S$ is thus a homomorphism of the Lie group $SO(3)$. As $F$ is a bijection, and $SO(3)$ has {\color{black} a trivial} center, $S$ is an isomorphism (automorphism) of the Lie group $SO(3)$. 
As in Section \ref{chap:Liegroups-with-biinvariant-metric},  denote by $\mu^+$ the (unique, up to a constant factor) biinvariant metric on $SO(3)$ with value $\mu_\epsilon^+:=\mu$ at $\epsilon$ and $\Theta:\mathfrak{so}(3) \to \mathfrak{so}^*(3), $ $x\mapsto \Theta(x)$, where $\Theta(x)=i_x\mu:=\mu(x,\cdot)$.  
We now consider the map 
$T: SO(3)\ltimes_{Ad^*}\mathfrak{so}^*(3)\to SE(3)$, $T(\sigma,f)=\Big(S(\sigma), F^{-1}( \Theta^{-1} (f))\Big)$. Clearly, $T$ is a bijection  and is further a homomorphism.
 Indeed, on the one hand, we have $T(\sigma,f)\;T(\tau,g)$ 
$ 
= \Big(S(\sigma)S(\tau), F^{-1}( \Theta^{-1} (f)) + S(\sigma)  F^{-1}( \Theta^{-1} (g)) \Big)$.  
 On the other hand, we have $T\Big((\sigma,f)(\tau,g)\Big)= \Big(S(\sigma\tau), F^{-1}( \Theta^{-1} (f)) + F^{-1}( \Theta^{-1} (Ad_\sigma^*g)) \Big)$. 
The equality $\Theta^{-1} \!(\!Ad_\sigma^*g)\!)\!\!=\!\!Ad_\sigma\Theta^{-1}\! (g)$ implies 
$F^{-1}\!(\!\Theta^{-1}\! (Ad_\sigma^*g)\!)\!\!=\!\! F^{-1}\!(\sigma F(F^{-1}\!(\Theta^{-1}\! (g)\!)\!)\sigma^{-1})$, 
 which reads $F^{-1}( \Theta^{-1} (Ad_\sigma^*g))= S(\sigma)(F^{-1}(\Theta^{-1} (g))) $.
 Altogether, we have 
\beqn T\Big((\sigma,f)(\tau,g)\Big)&=&\Big(S(\sigma\tau), F^{-1}( \Theta^{-1} (f)) + S(\sigma)(F^{-1}(\Theta^{-1} (g))) \Big)\nonumber\\
&=&T(\sigma,f)\;T(\tau,g).
\eeqn
By construction, $T$ is smooth, as $S$ is smooth and $F^{-1}$, $\Theta^{-1}$ are linear maps. Hence $T\circ \zeta: T^*SO(3)\to SE(3) $ is an isomorphism between the Lie groups  $T^*SO(3)$ and $ SE(3)$, where $\zeta: T^*SO(3)\to SO(3)\ltimes_{Ad^*} \mathfrak{so}^*(3)$, is  the right trivialization.
 Let us now show that $\bar \varphi:\hat{\mathbb H}_1\to SU(2)\ltimes \mathfrak{su}^*(2)$, $\bar\varphi (\hat Q) = \Big(\psi(Q), \Theta(\psi(\vec{u}))\Big)$,  {\color{black} where $\psi$ is as in (\ref{isomorphimSU(2)-quaternions})},  is a Lie group isomorphism between $\hat{\mathbb H}_1$ and $ SU(2)\ltimes \mathfrak{su}^*(2)$. {\color{black} For any
$\hat Q=Q+\mathfrak{e} \vec{u}Q,$ $\hat P=P+\mathfrak{e} \vec{v}P$ in $\hat{\mathbb H}_1,$} we have $\bar\varphi(\hat Q)\bar\varphi(\hat P)=
\Big(\psi(Q)\psi(P), \Theta(\psi(\vec{u})) + Ad_{\psi(Q)}^*\Theta(\psi(\vec{v}))\Big) $.
Applying the linearity of $\Theta$ and 
the identity $\Theta(Ad_\sigma x) = Ad_\sigma^*\Theta(x)$, the above come to
 $\bar \varphi(\hat Q)\bar\varphi(\hat P)=\Big(\psi(Q)\psi(P), \Theta(\psi(\vec{u}) + Ad_{\psi(Q)}\psi(\vec{v}))\Big)$. We further plug in the equality $\psi(\mathcal  R_Q\vec{v})=\psi(Q\vec{v}Q^{-1})= \psi(Q)\psi(\vec{v}) (\psi(Q))^{-1} = Ad_{\psi(Q)}\psi(\vec{v})$, to finally  get 
$\bar\varphi(\hat Q)\bar\varphi(\hat P) = \Big( \psi(QP), \Theta\psi(\vec{u}+\mathcal  R_Q\vec{v})\Big)=\bar\varphi(\hat Q\hat P) $. So $\bar\varphi$ is a Lie group homomorphism. By construction of $\Theta$ and $\psi$,   we see that $\bar\varphi$ is in fact, an isomorphism between the Lie groups $\hat{\mathbb H}_1$ and $ SU(2)\ltimes \mathfrak{su}^*(2)$. 
Again using the trivialization $\zeta_{SU}: T^*SU(2)\to SU(2)\ltimes_{Ad^*} \mathfrak{su}^*({\color{black}2})$,  we get the isomorphism $\varphi:=\zeta_{SU}^{-1}\circ \bar\varphi: \hat{\mathbb H}_1\to T^*SU(2) $  between the Lie groups  $\hat{\mathbb H}_1$ and $T^*SU(2)$.  More explicitly, for $ \hat Q=Q+\mathfrak{e}(\vec{u} Q) $ in  $ \hat{\mathbb H}_1$, we have $\varphi(\hat Q)= \zeta_{SU}^{-1}\Big(\psi(Q), \Theta(\psi(\vec{u})\Big)= \Big(\psi(Q), \Theta(\psi(\vec{u}))\circ T_{\psi(Q)}R_{\psi(Q)^{-1}}\Big)$. 
From Theorem \ref{thm:biinvariant-metric-on-dual}, the cotangent bundle $T^*SO(3)$  (resp. $T^*SU(2)$) and  the tangent bundle $TSO(3)$ (resp. $TSU(2)$)  of the compact Lie group $SO(3)$ (resp. $SU(2)$), are isomorphic. {\color{black} Note also that $ \hat{\mathbb H}_1\to  SU(2)\ltimes_{Ad} \mathfrak{su}(2),$ $Q+\mathfrak{e}\vec{u}\mapsto(\psi(Q),\psi(\vec{u}))$ is a Lie group isomorphism, from which one can also establish the above proof.}
\end{proof}
 Applying Theorem \ref{thm:biinvariant-metric-on-dual} to  $T^*SO(3)$ and further  identifying the latter with SE(3) according to Theorem \ref{theorem:commutingdiagramSO(3)}, we  recover some of  the results
on  SE$(3)$ obtained in \cite{pennec1}, \cite{Zefran}.  In particular, we recover  Theorem 3.6 of \cite{Zefran}.
\begin{theorem}\cite{Zefran} Let $\mu$  be a metric on SE($3$) such that every screw motion is a geodesic. Then $\mu$ is biinvariant and furthermore, its matrix is of the form (\ref{eq:biinvariantmetricTsl(2,1)}) in some basis of SE($3$). There is no Riemannian metric on SE($3$) for which every screw motion is a geodesic.
\end{theorem}
 Let us define 
$U:SO(2,1)\to SO(2,1)$, $U(\sigma):\mathbb R^3\to\mathbb R^3$,  where $ U(\sigma)x:=H^{-1}(\sigma H(x)\sigma^{-1}),$
with $H$ as in (\ref{defineH}),  (\ref{defineH(x)}).
It is easy to see that $U$ is an isomorphism (automorphism) of the Lie group $SO(2,1)$.
Let $\mu^+$ be the biinvariant Lorentzian metric on $SO(2,1)$,  $\mu_\epsilon^+:=\mu$ and  $\Theta:\mathfrak{so}(2,1) \to \mathfrak{so}^*(2,1), $ $x\mapsto \Theta(x):=\mu(x,\cdot)$.  
Further define 
$T': SO(2,1)\ltimes_{Ad^*}\mathfrak{so}^*(2,1)\to SE(2,1)$, $T'(\sigma,f)=\Big(U(\sigma), H^{-1}( \Theta^{-1} (f))\Big)$. Clearly, $T'$ is a bijection  and is further a homomorphism. Indeed, the product

$T'(\sigma,f)\;T'(\tau,g)$  reads
 $ \Big(U(\sigma)U(\tau), H^{-1}( \Theta^{-1} (f)) + U(\sigma)  H^{-1}( \Theta^{-1} (g)) \Big)$. On the other hand,  
$T'\Big((\sigma,f)(\tau,g)\Big)$ can be written as 
$$ 
T'\Big((\sigma,f)(\tau,g)\Big)=\Big(U(\sigma\tau), H^{-1}( \Theta^{-1} (f)) + H^{-1}( \Theta^{-1} (Ad_\sigma^*g)) \Big). 
$$
Applying the property 
$$\Theta^{-1} (Ad_\sigma^*g)=Ad_\sigma\Theta^{-1} (g)= \sigma\Theta^{-1} (g)\sigma^{-1}$$  to 
$T'\Big((\sigma,f)(\tau,g)\Big)$, the latter now becomes 
$$ \Big(U(\sigma\tau), H^{-1}( \Theta^{-1} (f)) + H^{-1}( \sigma H(H^{-1}(\Theta^{-1} (g)))\sigma^{-1}) \Big).$$ 
We rewrite it as
 $$
 T'\Big((\sigma,f)(\tau,g)\Big)=\Big(U(\sigma)U(\tau), H^{-1}( \Theta^{-1} (f)) + U(\sigma)(H^{-1}(\Theta^{-1} (g))) \Big), 
 $$ which is exactly $T'(\sigma,f)\;T'(\tau,g)$.
We get the isomorphism  $T'\circ \zeta$ between the Lie groups  $T^*SO(2,1)$ and $ SE(2,1) ,$ via the trivialization $\zeta: T^*SO(2,1)\to SO(2,1)\ltimes_{Ad^*} \mathfrak{so}^*(2,1).$  By {\color{black}Theorem \ref{thm:cotangentbundlerighttrivialization}}, the Lie groups $T^*SO(2,1)$  and $TSO(2,1)$ are isomorphic. {\color{black} We
construct  the Lie group isomorphism $\mathfrak{p}: \tilde{\mathbb S}_1\to  SL(2,\mathbb R)\ltimes_{Ad} \mathfrak{sl}(2,\mathbb R),$
$\mathfrak{p}(Q+\mathfrak{e}\vec{u}Q) = (\mathfrak{\omega}(Q), \mathfrak{\omega}(\vec{u}))$,
with $\mathfrak{\omega}$ as in (\ref{splitquaternions-as-2by2-matrices}). The trivialization $\xi: TSL(2,\mathbb R)\to  SL(2,\mathbb R)\ltimes_{Ad} \mathfrak{sl}(2,\mathbb R)$  gives the Lie group isomorphism $\xi^{-1}\circ \mathfrak{p}: \tilde{\mathbb S}_1\to  TSL(2,\mathbb R).$
 Since the Lie groups $TSL(2,\mathbb R)$ and $T^*SL(2,\mathbb R)$ are isomorphic (Theorem \ref{thm:cotangentbundlerighttrivialization}), we have proved the following. }
\begin{theorem}\label{theorem:cummutingdiagramSL(2)} The special affine Lie group $SE(2,1)=SO(2,1)\ltimes \mathbb R^3$  is isomorphic to both the tangent bundle $TSO(2,1)$ and  cotangent bundle $T^*SO(2,1)$ of $SO(2,1)$ endowed with their Lie group structure induced by the right trivialization.
The group $\tilde{\mathbb S}_1$ of unit dual split quaternions is isomorphic to {\color{black} both the tangent bundle $TSL(2,\mathbb R)$} and  cotangent bundle T$^*SL(2,\mathbb R)$ of $SL(2,\mathbb R),$ endowed with their Lie group structure induced by the right trivialization.
\end{theorem}  
On $\mathfrak{sl}(2,\mathbb R)$, consider the Lorentzian metric  $K_0(x,y) =4 Trace(xy)$, which is also its Killing form.
In the basis
$e_1=\frac{\sqrt{2}}{4}(E_{1,2}-E_{2,1})$, $e_2=\frac{\sqrt{2}}{4}(E_{1,1}-E_{2,2})$, $ e_3=\frac{\sqrt{2}}{4}(E_{1,2}+E_{2,1})$,
the matrix of $K_0$ is
diag($-1,1,1$)$=\mathbb I_{1,2}$. 
From Theorem \ref{thm:biinvariant-metric-on-dual}, 
every Cartan-Schouten metric on $ T^*SL(2,\mathbb R)$,  has the following nonzero coefficients, 
$\mu(e_2,e_2)=\mu(e_3,e_3)= -\mu(e_1,e_1)=s$, 
 $\mu(e_1,e_1^*)= \mu(e_2,e_2 ^*)=\mu(e_3,e_3^*)=t$,
 for some $s,t\in\mathbb R$, with $t\neq 0$,
 where   $(e_1^*,e_2^*,e_3^*)$ is  the dual basis of $(e_1,e_2,e_3)$.
Since the Lie algebras of  $T^*SL(2, \mathbb R)$ and $T^*SO(2, 1)$ are isomorphic, the above is also the expression of Cartan-Schouten metrics on $T^*SO(2, 1).$  More precisely, in the basis  
$e_1'=\frac{\sqrt{2}}{2}(E_{2,3}-E_{3,2})$, $e_2'=\frac{\sqrt{2}}{2}(E_{1,3}+E_{3,1})$, $ e_3'=-\frac{\sqrt{2}}{2}(E_{1,2}+E_{2,1})$, of  $\mathfrak{so}(2,1)$
 the matrix the Killing form is
diag($-1,1,1$)$=\mathbb I_{1,2}$. So in the basis $(e_j',(e_j')^*),$ j=1,2,3, every Cartan-Schouten metric on $ T^*SO(2,1)$, has the above form.
 Applying the isomorphism between $T^*SO(2, 1)$ and SE(2,1) proved in Theorem \ref{theorem:cummutingdiagramSL(2)}, we get the following.

\begin{theorem} Let $\mu$  be a metric on SE($2,1$) such that every screw motion on $\mathbb R^3_1$ is a geodesic. Then $\mu$ is biinvariant and furthermore, its matrix is of the form (\ref{eq:biinvariantmetricTsl(2,1)}) in some basis of SE($2,1$). There is no Riemannian metric on SE($2,1$) for which every screw motion is a geodesic.
\end{theorem}
\noindent
This could be seen as the  SE($2,1)$ version of the results on SE$(3)$ in \cite{pennec1}, \cite{Zefran}.

\begin{remark}Theorems  \ref{theorem:commutingdiagramSO(3)} and \ref{theorem:cummutingdiagramSL(2)} imply, in particular, that $SE(3)$, $SE(2,1)$,  $\hat{\mathbb H}_1$   and $\tilde{\mathbb S}_1$ all have exact symplectic structures and as Lie groups, they possess biinvariant metrics which all have  signature $(3,3).$
\end{remark}

\end{document}